\documentclass[12pt]{article}
\usepackage{enumerate}
\usepackage{amsfonts}
\usepackage{amsmath}
\usepackage{amssymb}
\usepackage{amsthm}

\def\D{\mathcal{D}}
\def\V{\mathcal{V}}
\def\N{\mathcal{N}}
\def\H{\mathcal{H}}
\def\K{\mathcal{K}}
\def\M{\mathcal{M}}

\def\S{\mathfrak{S}}
\def\F{\mathfrak{F}}
\def\C{\mathfrak{C}}
\def\T{\mathfrak{T}}
\def\B{\mathfrak{B}}

\newcommand{\rank}{\mathrm{rank}}

\newcommand{\id}{\mathrm{Id}}
\newcommand{\Tr}{\mathrm{Tr}}

\newcommand{\shs}{\hspace{1pt}}
\newcommand{\sn}{\|\hspace{-1pt}|}

\newcounter{defin}  \newcounter{lemma}  \newcounter{theorem}
\newcounter{property} \newcounter{corol}  \newcounter{remark} \newcounter{example}

\newenvironment{lemma}{\par\refstepcounter{lemma}
     \textbf{Lemma \thelemma.} }{\rm\par}
\newenvironment{theorem}{\par\refstepcounter{theorem}
     \textbf{Theorem \thetheorem.}\ }{\rm\par}
\newenvironment{property}{\par\refstepcounter{property}
     \textbf{Proposition \theproperty.}\ }{\rm\par}
\newenvironment{corollary}{\par\refstepcounter{corol}
     \textbf{Corollary \thecorol.} }{\rm\par}
\newenvironment{definition}{\par\refstepcounter{defin}
     \textbf{Definition \thedefin.}\ }{\rm\par}
\newenvironment{remark}{\par\refstepcounter{remark}
     \textbf{Remark \theremark.}}{\rm\par}

\begin{document}

\title{Operator \emph{E}-norms and their use}
\author{M.E. Shirokov\footnote{Steklov Mathematical Institute, RAS, Moscow, email:msh@mi.ras.ru}}
\date{}
\maketitle

\vspace{-10pt}

\begin{abstract}
We consider a family of equivalent norms (called operator \emph{E}-norms) on the algebra $\B(\H)$ of all bounded operators on a separable Hilbert space $\H$ induced by a positive  densely defined operator $G$ on $\H$. By choosing different generating operator $G$ one can obtain the operator \emph{E}-norms producing different topologies, in particular, the strong operator topology on bounded subsets of $\B(\H)$.

We obtain a generalised version of the Kretschmann-Schlingemann-Werner theorem, which shows continuity of the Stinespring representation of CP linear maps w.r.t. the energy-constrained $cb$-norm (diamond norm) on the set of CP linear maps and the operator \emph{E}-norm on the set of Stinespring operators.

The operator \emph{E}-norms induced by a positive operator $G$ are well defined for linear operators relatively bounded w.r.t. the operator $\sqrt{G}$ and the linear space of such operators equipped with any of these norms is a Banach space. We obtain explicit relations between the operator \emph{E}-norms and the standard characteristics of $\sqrt{G}$-bounded operators. The operator \emph{E}-norms allow to obtain simple upper bounds and continuity bounds for some functions depending on $\sqrt{G}$-bounded operators used in applications.
\end{abstract}

\tableofcontents

\section{Introduction}

The algebra $\B(\H)$ of all bounded linear operators on a separable Hilbert space $\H$, some its subalgebras and subspaces are basic objects in different
fields of modern mathematics and mathematical physics \cite{B&R,Paul,R&S}. In particular, $\B(\H)$ appears as an algebra of observables in the theory of quantum
systems while unital completely positive maps between such algebras called quantum channels play the role of dynamical maps in the Heisenberg picture \cite{H-SCI,Watrous,Wilde}.

The variety of different topologies on $\B(\H)$, relations between them and their "physical" sense are well known for anybody who is interested in functional analysis, theory of operator algebras, mathematical and theoretical physics.

In this article we describe families of norms on $\B(\H)$ producing different topologies on $\B(\H)$, in particular, the strong operator topology on bounded subsets of $\B(\H)$. These norms depending on a positive  densely defined operator ${G}$ and a positive parameter $E$ were introduced in \cite{CSR} for quantitative analysis of continuity of the Stinespring representation of a quantum channel with respect to the strong convergence of quantum channels and the strong operator convergence of Stinespring isometries.\footnote{Other applications of the operator \emph{E}-norms  are presented in the recent papers \cite{SPM,QDS,EPM}.}

Now we consider these norms (called the operator \emph{E}-norms) in more general context (assuming that ${G}$ is an arbitrary positive  operator). In Section 3 we consider equivalent definitions and basic properties of the operator \emph{E}-norms. We obtain explicit relations between the operator \emph{E}-norms and the equivalent norm on $\B(\H)$ also induced by a positive operator $G$ (which is commonly used in analysis of $\sqrt{G}$-bounded operators).

The operator \emph{E}-norms make it possible to obtain a generalization the Kretschmann-Schlingemann-Werner theorem. The original version of this theorem presented in \cite{Kr&W} shows continuity of the Stinespring representation of a completely positive (CP) linear map  with respect to the norm of complete boundedness ($cb$-norm in what follows)\footnote{It is also called the diamond norm in the quantum information theory \cite{Kit,Wilde}.} on the set of CP maps and the operator norm on the set of Stinespring operators. Our aim was to obtain a version of this theorem for other (weaker) topologies
on the sets of CP maps and corresponding Stinespring operators, in particular, for the strong convergence topology on the set of CP maps and the strong operator topology on the set of Stinespring operators. By using the operator \emph{E}-norms one can upgrade the proof of the Kretschmann-Schlingemann-Werner theorem without essential changes. The generalised version of this theorem and its corollaries are presented in Section 4.

In Section 5 the operator \emph{E}-norms induced by a positive operator $G$ are extended to linear operators  relatively bounded w.r.t. the operator $\sqrt{G}$. We prove that the linear space of such  operators equipped with any of these norms is a Banach space. Its subspace  consisting of all operators with zero $\sqrt{G}$-bound is the completion of $\B(\H)$ w.r.t. any of the operator \emph{E}-norms. We obtain explicit relations between the operator\break \emph{E}-norms and the standard characteristics of $\sqrt{G}$-bounded operators. \pagebreak

The operator \emph{E}-norms allow to obtain simple upper estimates and continuity bounds for some functions depending on $\sqrt{G}$-bounded operators used in applications.

As a basic example we consider the operators associated with the Heisenberg
Commutation Relation.

\section{Preliminaries}

Let $\mathcal{H}$ be a separable infinite-dimensional Hilbert space, $\mathfrak{B}(\mathcal{H})$
-- the algebra of all bounded operators on $\mathcal{H}$ with the operator norm $\|\!\cdot\!\|$ and $\mathfrak{T}(\mathcal{H})$ --
the Banach space of all trace-class operators on $\mathcal{H}$ with the trace norm $\|\!\cdot\!\|_1$ (the Schatten class of order 1) \cite{B&R,R&S}. Let
$\mathfrak{T}_{+}(\mathcal{H})$ be the cone of positive operators in
$\mathfrak{T}(\mathcal{H})$. Trace-class operators will be usually denoted by the Greek letters $\rho $, $\sigma $, $\omega $, ...
The closed convex subsets
$$
\T_{+,1}(\H)=\{\rho\in \T_+(\H)\,|\, \Tr\rho\leq 1\}\quad \textrm{and} \quad \S(\H)=\{\rho\in \T_+(\H)\,|\, \Tr\rho=1\}
$$
of the cone $\T_+(\H)$ are complete separable metric spaces with the metric defined by the trace norm.
Operators in $\S(\H)$ are called density operators or \emph{states}, since any $\rho$ in $\S(\H)$ determines a
normal state $A\mapsto \Tr A\rho$ on the algebra $\B(\H)$ \cite{B&R,H-SCI}. Extreme points of $\S(\H)$ are 1-rank projectors called \emph{pure states}.

Denote by $I_{\H}$ the unit operator on a Hilbert space
$\mathcal{H}$ and by $\id_{\mathcal{\H}}$ the identity
transformation of the Banach space $\mathfrak{T}(\mathcal{H})$.

We will use the Dirac notations $|\varphi\rangle$, $|\psi\rangle\langle\varphi|$,... for vectors and operators of rank 1 on a Hilbert space (in this notations the action of an operator $|\psi\rangle\langle\varphi|$ on a vector $|\chi\rangle$ gives the vector $\langle\varphi|\chi\rangle|\psi\rangle$) \cite{H-SCI}.

We will pay a special attention to the class of  unbounded densely defined positive operators on $\H$ having discrete spectrum of finite multiplicity.
In Dirac's notations any such operator ${G}$ can be represented as follows
\begin{equation}\label{H-rep}
G=\sum_{k=0}^{+\infty} E_k|\tau_k\rangle\langle \tau_k|
\end{equation}
on the domain $\mathcal{D}(G)=\{ \varphi\in\H\,|\,\sum_{k=0}^{+\infty} E^2_k|\langle\tau_k|\varphi\rangle|^2<+\infty\}$, where
$\left\{\tau_k\right\}_{k=0}^{+\infty}$ is the orthonormal basis of eigenvectors of ${G}$
corresponding to the nondecreasing sequence $\left\{\smash{E_k}\right\}_{k=0}^{+\infty}$ of eigenvalues
tending to $+\infty$. We will use the following (cf.\cite{W-EBN})\smallskip

\begin{definition}\label{D-H}
An operator ${G}$ having representation (\ref{H-rep}) is called \emph{discrete}.
\end{definition}\smallskip

The set $\S(\H)$ is compact if and only if $\dim \H<+\infty$. We will  use the following\smallskip

\begin{lemma}\label{Comp} \cite{H-c-w-c} \emph{If ${G}$ is a discrete unbounded operator on $\H$ then the set of states $\rho$ in $\S(\H)$ satisfying the inequality $\Tr {G}\rho\leq E$ is compact for any $E\geq\inf_{\|\varphi\|=1}\langle\varphi|{G}|\varphi\rangle$.}
\end{lemma}\smallskip

We will also use the following simple lemma.

\begin{lemma}\label{WL} \cite{W-CB} \emph{If $f$ is a concave nonnegative function on $[0,+\infty)$ then for any positive $x< y$ and any $z\geq0$ the inequality $\,xf(z/x)\leq yf(z/y)\,$ holds.}
\end{lemma}

\section{Operator \emph{E}-norms on $\B(\H)$}

Let ${G}$ be a positive semidefinite operator on $\H$ with a dense domain $\mathcal{D}({G})$ such that
\begin{equation}\label{H-cond}
\inf\left\{\shs\|G\varphi\|\,|\,\varphi\in\mathcal{D}({G}),\|\varphi\|=1\shs\right\}=0.
\end{equation}
We will assume that for any positive operator $\rho$ in $\T(\H)$ the value of $\Tr {G}\rho$ (finite or infinite) is defined as $\sup_n\Tr P_nG\rho$, where $P_n$ is the spectral projector of $G$ corresponding to the interval $[0,n]$.\smallskip

For given $E>0$ consider the function on $\B(\H)$ defined as
\begin{equation}\label{ec-on-b}
 \|A\|^G_E\doteq \sup_{\varphi\in\H_1, \langle\varphi|G|\varphi\rangle\leq E}\|A\varphi\|,
\end{equation}
where $\H_1$ is the unit sphere in $\H$ and it is assumed that $\,\langle\varphi|G|\varphi\rangle=\|\sqrt{G}\varphi\|^2\,$ if $\varphi$ lies in $\D(\sqrt{G})$ and $\,\langle\varphi|G|\varphi\rangle=+\infty\,$ otherwise. This function can be also defined as
\begin{equation}\label{ec-on}
 \|A\|^{G}_E\doteq \sup_{\substack{\rho\in\mathfrak{S}(\mathcal{H}):
\Tr {G}\rho\leq E}}\sqrt{\Tr A\rho A^*},
\end{equation}
where the supremum is over all states $\rho$ in $\S(\H)$ satisfying the inequality $\Tr {G}\rho\leq E$.\footnote{In the previous versions of this posting the coincidence of the r.h.s. of (\ref{ec-on-b}) and (\ref{ec-on}) was conjectured, but it was proved only under the assumption that
the operator $G$ is discrete (Definition \ref{D-H}).}

The coincidence of the r.h.s. of (\ref{ec-on-b}) and (\ref{ec-on}) for any $A\in\B(\H)$ is shown in \cite{W&Sh}.

It is easy to see that the function $A\mapsto\|A\|^{G}_E$  is a  norm on $\B(\H)$.
Definition (\ref{ec-on-b}) shows the sense of the norm $\|\cdot\|^G_E$ (as a constrained version of the operator norm $\|\cdot\|$) while definition
(\ref{ec-on}) is more convenient for studying its analytical properties. In particular, by using definition
(\ref{ec-on}) the following proposition is proved in \cite{CSR}.\footnote{In \cite{CSR} condition (\ref{H-cond}) was not assumed. We use this assumption here, since it simplifies analysis of the norms $\|\!\cdot\!\|^{G}_E$ without reduction of generality (note that $\|A\|^{{G}+\lambda I}_E=\|A\|^{G}_{E-\lambda}$ for all $A$ and $\lambda>0$).} \smallskip

\begin{property}\label{ec-on-p1}
\emph{For any  operator  $\,A\in\B(\H)$ the following properties hold:}
\begin{enumerate}[a)]
  \item \emph{$\|A\|^{G}_E$ tends to $\|A\|$ as $E\rightarrow+\infty$;}
  \item \emph{the function $E\mapsto\left[\|A\|^{G}_E\right]^p$ is concave and nondecreasing on $\,\mathbb{R}_+$ for $p\in(0,2]$;}
  \item \emph{$\|A\varphi\|\leq K_{\varphi}\|A\|^{G}_E$ for any unit vector $\varphi$ in $\D(\sqrt{G})$, where $K_{\varphi}=\max\{1, \|\sqrt{G}\varphi\|/\sqrt{E}\}$.}
\end{enumerate}
\end{property}

We will call the norms $\|\!\cdot\!\|^{G}_E$ the \emph{ operator E-norms} on $\B(\H)$. Property b) in Proposition \ref{ec-on-p1} shows that
\begin{equation}\label{E-n-eq}
\|A\|^{G}_{E_1}\leq \|A\|^{G}_{E_2}\leq \sqrt{E_2/E_1}\|A\|^{G}_{E_1}\quad\textrm{ for any } E_2>E_1>0.
\end{equation}
Hence for given operator ${G}$ all the norms $\|\!\cdot\!\|^{G}_{E}$, $E>0$, are equivalent on $\B(\H)$.

\smallskip

\begin{remark}\label{en-gen}
The definition of the operator \emph{E}-norm is obviously generalized to operators between different Hilbert spaces $\H$ and $\K$. It is easy to see that all the above and below results concerning properties of the operator \emph{E}-norms
remain valid (with obvious modifications) for this generalization. $\square$
\end{remark}\smallskip

Since the set $\D(\sqrt{G})$ is dense in $\H$, property c) in Proposition \ref{ec-on-p1} shows that the topology generated by any of the norms $\|\!\cdot\!\|^{G}_{E}$ on bounded subsets of $\B(\H)$ is not weaker than the strong operator topology. On the other hand, it is not stronger than the norm topology on $\B(\H)$.
The following proposition characterizes these extreme cases.
\smallskip

\begin{property}\label{ec-on-p2}
A) \emph{The norm $\,\|\!\cdot\!\|^{G}_E$, $E>0$, is  equivalent to the operator norm $\,\|\!\cdot\!\|$ on $\B(\H)$ if and only if the operator ${G}$ is bounded.} \smallskip

B) \emph{The norm $\,\|\!\cdot\!\|^{G}_E$, $E>0$, generates the strong operator topology on bounded subsets of $\B(\H)$
if and only if $\,{G}$ is an unbounded discrete operator (Definition \ref{D-H}).}
\end{property}\smallskip

\emph{Proof.} A) If ${G}$ is a bounded operator then $\,\|\!\cdot\!\|^{G}_E=\|\!\cdot\!\|$ for any $E\geq \|G\|$.

If ${G}$ is a unbounded operator and $P_n$ is the spectral projector of ${G}$ corresponding to the interval $[n,+\infty)$ then $\|P_n\|=1$ for all $n$. By noting that $\Tr P_n\rho\leq E/n$ for any state $\rho$ such that $\Tr {G}\rho\leq E$, it is easy to see that $\,\|P_n\|^{G}_E\rightarrow 0$  as $n\rightarrow+\infty$. \smallskip

B) The "if" part of this assertion is proved in
\cite{CSR}.

Assume there is a spectral projector of the operator ${G}$ corresponding to a finite interval $[0, E_0]$ with infinite-dimensional range $\H_0$.
Since $\langle\varphi|G|\varphi\rangle\leq E_0$ for any unit vector $\varphi$ in $\H_0$, we have $\|A\|^{G}_E=\|A\|$ for any $A\in\B(\H_0)$ and $E>E_0$. So, any of the norms $\,\|\!\cdot\!\|^{G}_E$, $E>0$, generates the norm topology on $\B(\H_0)$ in this case $\square$. \smallskip

Different types of operator convergence can be obtained by using the norm $\|\!\cdot\!\|_E^{G}$ induced by  different operators ${G}$. \smallskip

\textbf{Example.} Let $\H=\H_1\oplus\H_2$ and ${G}={G}_1\oplus {G}_2$, where ${G}_k$ is a  positive  densely defined operator on a separable Hilbert space $\H_k$ satisfying condition (\ref{H-cond}), $k=1,2$. By using definition (\ref{ec-on}) and the triangle inequality it is easy to show that
\begin{equation}\label{E-ineq}
  \sqrt{p\!\left[\|AP_1\|^{{G}_1}_E\right]^2+(1-p)\!\left[\|AP_2\|^{{G}_2}_E\right]^2}\leq \|A\|^{G}_E\leq \|AP_1\|^{{G}_1}_E+\|AP_2\|^{{G}_2}_E
\end{equation}
for any $p\in[0,1]$, where $P_k$ is the projector on the subspace $\H_k$ and $\|AP_k\|^{{G}_k}_E$, $k=1,2$, are defined in accordance with Remark \ref{en-gen}.

Assume that ${G}_1$ is a discrete unbounded operator (Def.\ref{D-H}) and ${G}_2$ is a bounded operator. Then it follows from (\ref{E-ineq}) and Proposition \ref{ec-on-p2} that
$$
\left\{\|\!\cdot\!\|^{G}_E\shs\textrm{-}\lim_{n\rightarrow\infty} A_n=A_0\right\}\quad \Leftrightarrow\quad \left\{{s.o.}\shs\textrm{-}\lim_{n\rightarrow\infty} A_nP_1=A_0P_1\right\} \wedge \left\{\|\!\cdot\!\|\shs\textrm{-}\lim_{n\rightarrow\infty} A_nP_2=A_0P_2\right\}
$$
for a bounded sequence $\{A_n\}\subset\B(\H)$, where ${s.o.}\shs\textrm{-}\lim$ denotes the limit  w.r.t the strong operator topology. So, in this case the norm $\|\!\cdot\!\|^{G}_E$ generates a "hybrid" topology on bounded subsets of
$\B(\H)$ -- some kind of the Cartesian product of the strong operator and the norm topologies.

\subsection{Equivalent definitions and equivalent norms}

Recall that $\T_{+,1}(\H)$ denotes the positive part of the unit ball in $\T(\H)$. Denote by $\H_{\leq 1}$  the unit ball in $\H$.\smallskip

\begin{property}\label{en-def}  A) \emph{For any $A\in\B(\H)$ and $E>0$ the following expressions hold
\begin{equation}\label{ec-on+}
 \|A\|^{G}_E=\sup_{\substack{\varphi\in\H_{\leq 1}:\langle\varphi|G|\varphi\rangle\leq E}}\|A\varphi\|=\sup_{\substack{\rho\in\T_{+,1}(\H):
\Tr {G}\rho\leq E}}\sqrt{\Tr A\rho A^*},\qquad \forall A\in \B(\H),
\end{equation}
i.e. the suprema in definitions (\ref{ec-on-b}) and (\ref{ec-on}) can be taken, respectively, over all vectors in $\H_{\leq 1}$ satisfying the condition
$\,\langle\varphi|G|\varphi\rangle\leq E$ and over all operators in $\T_{+,1}(\H)$ satisfying the condition $\,\Tr {G}\rho\leq E$.} \smallskip

\noindent B) \emph{If the operator $G$ is unbounded then for any $A\in\B(\H)$ and $E>0$  the conditions $\,\langle\varphi|G|\varphi\rangle\leq E$ in (\ref{ec-on-b}) and $\,\Tr {G}\rho\leq E$ in (\ref{ec-on}) and can be replaced, respectively, by the conditions $\,\langle\varphi|G|\varphi\rangle=E$ and $\,\Tr G\rho=E$.} \smallskip

\noindent C) \emph{If $G$ is a discrete unbounded operator (Def.\ref{D-H}) then for any $A\in\B(\H)$ and $E>0$ the suprema in (\ref{ec-on-b}) and (\ref{ec-on}) are attainable. Moreover, if $\|A\|^{G}_{E}<\|A\|$  then the suprema in (\ref{ec-on-b}) and (\ref{ec-on}) are attained, respectively, at unit vector $\varphi_0$ in $\H$ such that $\,\langle\varphi_0|G|\varphi_0\rangle=E$ and at a state
$\rho_0$ in $\S(\H)$ such that $\,\Tr G\rho_0=E$.}
\end{property}\smallskip

\begin{remark}\label{en-def-r} It is easy to see that assertion A of Proposition \ref{en-def} is not valid if the operator $G$ doesn't satisfy condition (\ref{H-cond}).
\end{remark}\smallskip

\emph{Proof of Proposition \ref{en-def}.} A) It suffices to show that
the last expression in (\ref{ec-on+}) does not exceed $\|A\|^{G}_E$. Let $\rho$ be an operator in $\T_{+,1}(\H)$ such that $\,\Tr {G}\rho\leq E$ and $r=\Tr\rho$. Then $\hat{\rho}\doteq r^{-1}\rho$ is a state such that
$\,\Tr {G}\hat{\rho}\leq E/r$. So, by using concavity of the function $E\rightarrow\left[\|A\|^{G}_E\right]^2$ and Lemma \ref{WL} in Section 2 we obtain
$$
\Tr A\rho A^*=r\Tr A\hat{\rho} A^*\leq r\left[\|A\|^{G}_{E/r}\right]^2\leq \left[\|A\|^{G}_{E}\right]^2.
$$

B)\footnote{If $\|A\|^{G}_E<\|A\|$ then this assertion can be derived from properties a) and b) in Proposition \ref{ec-on-p1}.} Show first that the inequality $\Tr {G}\rho\leq E$ can be replaced by the equality $\Tr {G}\rho=E$ in (\ref{ec-on}). Assume that there exist  $E>0$ and $A\in\B(\H)$ such that
\begin{equation}\label{assum}
\sup_{\substack{\rho\in\mathfrak{S}(\mathcal{H}):
\Tr {G}\rho=E}}\sqrt{\Tr A\rho A^*}\leq \|A\|^{G}_E-\varepsilon
\end{equation}
for some $\varepsilon>0$. Let $\rho_\varepsilon$ be a state such that
$\sqrt{\Tr A\rho_\varepsilon A^*}>\|A\|^{G}_E-\varepsilon/2$ and $\Tr G\rho_\varepsilon<E$. For each natural $n>E$ there exist a state $\sigma_n$ and a number $p_n\in(0,1)$ such that $\Tr G\sigma_n\in(n,+\infty)$ and $\Tr G\varrho_n=E$, where $\varrho_n=(1-p_n)\rho_\varepsilon+p_n\sigma_n$. It is clear that $p_n\rightarrow 0$ as $n\rightarrow+\infty$. Hence $\Tr A\varrho_n A^*$ tends to $\Tr A\rho_\varepsilon A^*$ contradicting  (\ref{assum}).

For any $\varepsilon>0$ let $\rho_\varepsilon$ be a state such that
$\Tr A\rho_\varepsilon A^*>[\|A\|^{G}_E-\varepsilon]^2$ and $\Tr G\rho_\varepsilon=E$. By Corollary 1 in \cite{W&Sh} there is a probability
measure $\mu$ on $\S(\H)$ supported by pure states  such that
\begin{equation}\label{p-s-d}
\rho_\varepsilon=\int \sigma\mu(d\sigma)\quad \textrm{and}\quad \Tr H\sigma=E\;\;\textrm{for}\;\;\mu\textrm{-almost all}\;\;\sigma.
\end{equation}
Since that function $\sigma\mapsto\Tr A\sigma A^*$ is affine and continuous, we have
$$
\int \Tr A\sigma A^*\mu(d\sigma)=\Tr A\rho_\varepsilon A^*>[\|A\|^{G}_E-\varepsilon]^2.
$$
It shows existence of a pure state $\sigma_\varepsilon$ such that $\Tr A\sigma_\varepsilon A^*>[\|A\|^{G}_E-\varepsilon]^2$ and $\Tr H\sigma_\varepsilon=E$.\smallskip

C) In this case the set of pure states $\rho$ satisfying the condition $\Tr {G}\rho\leq E$ is compact by Lemma \ref{Comp} in Section 2. Hence, by the coincidence of the r.h.s. of (\ref{ec-on-b}) and (\ref{ec-on}), the supremum in (\ref{ec-on}) is attained at some pure state $\rho_0$. By Proposition \ref{ec-on-p1} the condition $\|A\|^{G}_E<\|A\|$ shows that $\|A\|^{G}_{E'}<\|A\|^{G}_E$ for any $E'<E$. Hence $\Tr G\rho_0=E$.  $\square$
\medskip

Consider the following norm on $\B(\H)$ depending on a positive parameter $E$:
\begin{equation}\label{eq-norms-2}
\sn A\sn^{G}_{E}=\sup\left\{\|A\varphi\|\left|\, \varphi\in \D(\sqrt{G}),\, \|\varphi\|^2+\langle\varphi|G|\varphi\rangle/E\leq 1\right.\right\}.
\end{equation}
This norm naturally appears in analysis of operators relatively bounded w.r.t. the operator $\sqrt{G}$ (see Section 5).  We will obtain relations between the norms
$\|\cdot\|^{G}_{E}$ and $\sn\cdot\sn^{G}_{E}$ assuming that $G$ is an arbitrary positive operator satisfying condition (\ref{H-cond}).\smallskip

By using definitions (\ref{ec-on-b}) and (\ref{eq-norms-2}) it is easy to show that
\begin{equation}\label{one}
\sqrt{1/2}\|A\|^{G}_{E}\leq \sn A\sn^{G}_{E}\leq\|A\|^{G}_{E}
\end{equation}
for any $E>0$. This inequality and inequality (\ref{E-n-eq}) imply that
\begin{equation}\label{E-n-eq+}
\sn A\sn^{G}_{E_1}\leq \sn A\sn^{G}_{E_2}\leq \sqrt{2E_2/E_1}\sn A\sn^{G}_{E_1}\quad\textrm{ for any } E_2>E_1>0
\end{equation}
and any $A\in\B(\H)$. The above inequalities show that all the norms in the families $\{\|\cdot\|^{G}_{E}\}_{E>0}$ and $\{\sn \cdot\sn^{G}_{E}\}_{E>0}$ are equivalent to each other. \smallskip

In fact, all the norms $\|\cdot\|^{G}_{E}$ and $\sn \cdot\sn^{G}_{E}$ are equivalent on the space of all linear operators on $\H$ relatively bounded w.r.t. the operator $\sqrt{G}$, and each of these norms makes this space a Banach space. Moreover,
the functions $E\mapsto\|A\|^G_E$ and $E\mapsto\sn A\sn^G_E$ are completely determined by each other (Remark \ref{comp-d-r+} in Section 5).
The main advantages of the operator \emph{E}-norm $\|\cdot\|^{G}_{E}$ in comparison with the norm  $\sn \cdot\sn^{G}_{E}$ are the following:
\begin{itemize}
  \item the concavity of the function $E\mapsto\left[\|A\|^{G}_{E}\right]^p$ for any $p\in(0,2]$;\footnote{The function $E\mapsto\left[\sn A\sn^{G}_{E}\right]^p$ is not concave in general for any $p\in(0,2]$ by Remark \ref{nonconc} below;}
  \item the appearance of the norm  $\|\cdot\|^{G}_{E}$ in the generalized Kretschmann-Schlingemann-Werner theorem (Section 4);
  \item the simple estimation of $\|\Phi(A)\|^{G}_E$ via $\|A\|^{G}_E$, where $\Phi:\B(\H)\rightarrow\B(\H)$ is any 2-positive linear map satisfying the conditions of Proposition \ref{bp-en-0}E (Section 3.2).
\end{itemize}

\begin{remark}\label{nonconc} To show that the function $E\mapsto\left[\sn A\sn^{G}_{E}\right]^p$ is not concave in general for any $p\in(0,2]$ it suffices to consider two-dimensional Hilbert space $\H=\mathbb{C}^2$ and the operators
      $$
G=\left[\begin{array}{ll}
        1 & 0\\
        0 & 0
        \end{array}\right]\quad \textrm{and}   \quad A=\left[\begin{array}{ll}
        \sqrt{2} & 0\\
        0 & 1
        \end{array}\right].
$$
It is easy to see that $\,\sn A\sn^{G}_{E}=1\,$ if $\,E\in(0,1]\,$ and $\,\sn A\sn^{G}_{E}=\sqrt{2E/(E+1)}\,$ if $\,E>1$.
\end{remark}\smallskip

\subsection{Basic properties of the operator \emph{E}-norms}

In the following proposition we collect properties of the operator \emph{E}-norms used below.\smallskip

\begin{property}\label{bp-en-0} \emph{Let ${G}$ be a positive  densely defined operator  on a Hilbert space $\H$
satisfying condition (\ref{H-cond}) and $E>0$.}\medskip

\noindent A) \emph{$\|A\|^{G}_E=\||A|\|^{G}_E\leq \sqrt{\|A^*A\|^{G}_E}$ for all $A\in\B(\H)$  but $\,\|A^*\|^{G}_E\neq\|A\|^{G}_E$ in general;}\smallskip

\noindent B) \emph{For arbitrary operators $A$ and $B$ in $\B(\H)$ the following inequalities hold
$$
m(A)\|B\|^{G}_E\leq \|AB\|^{G}_E\leq \|A\|\|B\|^{G}_E,
$$
where $m(A)$ is the infimum of the spectrum of the operator $|A|=\sqrt{A^*A}$.}\smallskip

\noindent C) \emph{For arbitrary operators $A$ and $B$ in $\B(\H)$ such that $\,\langle A\varphi |B\varphi\rangle=0$ for any $\varphi\in\H$ the following inequalities hold}
$$
\max\!\shs\left\{\|A\|^{G}_{E},\|B\|^{G}_{E}\shs\right\}\leq \|A+B\|^{G}_{E}\leq \sqrt{\left[\|A\|^{G}_{E}\right]^2+\left[\|B\|^{G}_{E}\right]^2}.
$$

\noindent D) \emph{For any operator $\rho$ in $\T_{+,1}(\H)$ with finite $E_{\rho}\doteq\Tr {G}\rho$ and arbitrary operators $A$ and $B$ in $\B(\H)$ the following inequalities hold}
$$
|\Tr A\rho B^*|\leq \|A\rho B^*\|_1 \leq \|A\|^{G}_{E_{\rho}}\|B\|^{G}_{E_{\rho}}.
$$

\noindent E) \emph{For any 2-positive map $\Phi:\B(\H)\rightarrow\B(\H)$ such that $\,\Phi(I_{\H})\leq I_{\H}\,$ having the predual map\footnote{The map $\Phi_*$ is  defined by the relation $\Tr \Phi(A)\rho=\Tr A\Phi_*(\rho)$ for all $A\in\B(\H)$ and $\rho\in\T(\H)$. Existence of $\Phi_*$ is equivalent to normality of the map $\Phi$, which means that $\Phi(\sup_{\lambda}A_\lambda)=\sup_{\lambda}\Phi(A_{\lambda})$ for any increasing net $A_\lambda$ of positive operators in $\B(\H)$ \cite{B&R}.} $\Phi_*:\T(\H)\rightarrow\T(\H)$ with finite $\,Y_{\Phi}(E)\doteq \sup\!\left\{\Tr {G}\Phi_*(\rho)\,|\,\rho\in\mathfrak{S}(\mathcal{H}),\Tr {G}\rho\leq E\shs\right\}$ and arbitrary operator $A$  in $\B(\H)$ the following inequalities hold }\footnote{If ${G}$ is a Hamiltonian of a quantum system described by the space $\H$ and $\Phi$ is a quantum channel (in the Heisenberg picture) then $Y_{\Phi}(E)/E$ is the energy amplification factor of $\Phi$.}
$$
\|\Phi(A)\|^{G}_E\leq\sqrt{\|\Phi(I_{\H})\|}\,\|A\|^{G}_{Y_{\Phi}(E)}\leq \sqrt{\|\Phi(I_{\H})\|K_{\Phi}}\,\|A\|^{G}_{E},\qquad  K_{\Phi}=\max\{1, Y_{\Phi}(E)/E\}.
$$
\end{property}\smallskip

Proposition \ref{bp-en-0} shows that the linear transformations  $A\mapsto BA$ and $A\mapsto\Phi(A)$ of $\B(\H)$, where $B\in\B(\H)$ and $\Phi$ is a map with the properties pointed in part E, are bounded operators w.r.t. the norm $\|\!\cdot\!\|^{G}_E$ (in contrast to the transformation $A\mapsto AB$).

\smallskip

\emph{Proof.} A) The equality $\|A\|^{G}_E=\||A|\|^{G}_E$ is obvious. The inequality $\|A\|^{G}_E\leq\sqrt{\|A^*A\|^{G}_E}$ follows from the operator Cauchy-Schwarz inequality
$$
[\Tr A^*A\rho]^2\leq [\Tr [A^*A]^2\rho][\Tr\rho].
$$
To show that $\|A^*\|^{G}_E$ may not coincide with $\|A\|^{G}_E$ take any operator ${G}$ having form (\ref{H-rep}). It is easy to see that
$\||\tau_0\rangle\langle\tau_n|\|^{G}_E=\sqrt{E/E_n}\,$ while $\||\tau_n\rangle\langle\tau_0|\|^{G}_E=1\,$ for all $E>0$.
\smallskip

B)  This assertion follows directly from definition (\ref{ec-on}) of the operator \emph{E}-norm.\smallskip

C) This assertion follows directly from definition (\ref{ec-on-b}) of the operator \emph{E}-norm.\smallskip


D) The first inequality is obvious. Let  $U$ be the partial isometry from the polar decomposition of $A\rho B^*$, i.e. $A\rho B^*=U|A\rho B^*|$. By using the operator Cauchy-Schwarz inequality we obtain
$$
\|\Tr A\rho B^*\|^2_1=[\Tr U^*\!A\rho B^*]^2\leq [\Tr UU^*\!A\rho A^*][\Tr B\rho B^*]\leq [\Tr A\rho A^*][\Tr B\rho B^*],
$$
where that last inequality is due to the fact that $UU^*\leq I_{\H}$. By Proposition \ref{en-def}A the right hand side of this inequality does not exceed $\left[\|A\|^{G}_{E_{\rho}}\|B\|^{G}_{E_{\rho}}\right]^2$.\smallskip

E) By Kadison's inequality and Proposition \ref{en-def}A we have
$$
\Tr [\Phi(A)]^*\Phi(A)\rho\leq\|\Phi(I_{\H})\|\Tr\Phi(A^*A)\rho=\|\Phi(I_{\H})\|\Tr A^*A\Phi_*(\rho)\leq \|\Phi(I_{\H})\| \left[\|A\|^{G}_{Y_{\Phi}(E)}\right]^2
$$
for any $A\in\B(\H)$ and any $\rho\in\S(\H)$ such that $\Tr {G}\rho\leq E$ (since the condition $\,\Phi(I_{\H})\leq I_{\H}\,$ guarantees  that $\Phi_*(\rho)\in\T_{+,1}(\H)$). This implies the first inequality. The second inequality follows from (\ref{E-n-eq}). $\square$

\subsection{Properties of the \emph{E}-norms related to tensor products}

If ${G}_1$ and ${G}_2$ are positive  densely defined operators on Hilbert spaces $\H_1$ and $\H_2$
satisfying condition (\ref{H-cond}) then
${G}_{12}={G}_{1}\otimes I_{2}+ I_{1}\otimes {G}_{2}$ is an operator on the Hilbert space $\H_{12}=\H_1\otimes\H_2$ with the same properties.\footnote{If ${G}_1$ and ${G}_2$ are Hamiltonians of quantum systems $1$ and $2$ described by the spaces $\H_1$ and $\H_2$ then ${G}_{12}$ is the Hamiltonian of the composite quantum system $12$ \cite{H-SCI}.}\footnote{Here and in what follows we write $I_X$ instead of $I_{\H_X}\!$ (where $X=1,2,A,B,..$) to simplify notations.} The following proposition contains several estimates for the operator \emph{E}-norms of product operators used in Section 4.\smallskip

\begin{property}\label{pen-tp}  \emph{Let ${G}_1$, ${G}_2$ and ${G}_{12}$ be the operators described above.}\smallskip

\noindent A) \emph{For arbitrary operator $A$ in $\B(\H_1)$ the following equalities hold}
$$
\|A\otimes I_{2}\|^{{G}_{12}}_E=\|A\otimes I_{2}\|^{{G}_{1}\otimes \shs I_{2}}_E=\|A\|^{{G}_1}_E
$$

\noindent B) \emph{For arbitrary operators $A\in\B(\H_1)$ and $B\in\B(\H_2)$ the following inequalities hold}
\begin{equation}\label{pen-tp-1}
\sup_{x\in(0,E)}\|A\|^{{G}_1}_{x}\|B\|^{{G}_2}_{E-x}\leq \|A\otimes B\|^{{G}_{12}}_{E} \leq\sup_{x\in(0,E)} \sqrt{\|A^*A\|^{{G}_1}_{x}\|B^*B\|^{{G}_2}_{E-x}},
\end{equation}
\emph{and}
\begin{equation}\label{pen-tp-2}
\|A\otimes B\|^{{G}_{12}}_{E} \leq\min\left\{\|A\|^{{G}_1}_{E}\|B\|,\|A\|\|B\|^{{G}_2}_{E}\right\}.
\end{equation}
\end{property}\smallskip

Note that the lower and upper bounds in (\ref{pen-tp-1}) and the r.h.s. of (\ref{pen-tp-2}) tend to $\|A\|\|B\|=\|A\otimes B\|=\lim_{E\to+\infty}\|A\otimes B\|^{{G}_{12}}_{E}$ as $E\to+\infty$.\smallskip

\emph{Proof.} A) It suffices to note that $\Tr [|A|^2\otimes I_{2}]\rho_{12}=\Tr |A|^2\rho_{1}$ and
that $\Tr {G}_{12}\rho_{12}=\Tr {G}_{1}\rho_{1}+\Tr {G}_{2}\rho_{2}$ for any state $\rho_{12}\in\S(\H_{12})$, where $\rho_1=\Tr_{\H_2}\rho_{12}$ and $\rho_2=\Tr_{\H_1}\rho_{12}$ are the partial states of $\rho_{12}$.

B) For each $x\in(0,E)$ and any $\varepsilon>0$ there exist states $\rho_1$ in $\S(\H_1)$ and $\rho_2$ in $\S(\H_2)$
such that $\sqrt{\Tr |A|^2\rho_1}>\|A\|^{{G}_1}_x-\varepsilon$, $\Tr {G}_1\rho_1\leq x$, $\sqrt{\Tr |B|^2 \rho_2}>\|B\|^{{G}_2}_{E-x}-\varepsilon$ and  $\Tr {G}_2\rho_2\leq E-x$. Then
$\Tr {G}_{12}[\rho_1\otimes\rho_2]\leq x+E-x=E$ and  $\sqrt{\Tr [|A|^2\otimes|B|^2][\rho_1\otimes\rho_2]}\geq[\|A\|^{{G}_1}_x-\varepsilon][\|B\|^{{G}_2}_{E-x}-\varepsilon]$. Since $\varepsilon$ is arbitrary, this implies the left inequality in (\ref{pen-tp-1}).\smallskip

By the operator Cauchy-Schwarz inequality for any state $\rho_{12}$ in $\S(\H_{12})$  we have
$$
\begin{array}{c}
\Tr [|A|^2\otimes |B|^2]\rho_{12}\leq \sqrt{\Tr [|A|^4\otimes I_{2}]\rho_{12}}\sqrt{\Tr [I_{1}\otimes |B|^4]\rho_{12}}\\\\=
\sqrt{\Tr|A|^4\rho_{1}}\sqrt{\Tr |B|^4\rho_{2}}\leq \||A|^2\|^{{G}_1}_{\Tr {G}_1\rho_1}\||B|^2\|^{{G}_2}_{\Tr {G}_2\rho_2}.
\end{array}
$$
Since $\Tr {G}_{12}\rho_{12}=\Tr {G}_{1}\rho_{1}+\Tr {G}_{2}\rho_{2}$,  this implies the right inequality in (\ref{pen-tp-1}).\smallskip

To prove inequality (\ref{pen-tp-2}) it suffices to note that
$$
A\otimes B=[A\otimes I_{2}][I_{1}\otimes B]=[I_{1}\otimes B][A\otimes I_{2}]
$$
and to apply Proposition \ref{bp-en-0}B and part A of this proposition. $\square$

\section{The \emph{E}-version of the Kretschmann-Schlingemann-Werner theorem}

In this section we consider application of the operator $E$-norms to the theory of completely positive (CP) linear maps between Banach spaces of trace class
operators on  separable Hilbert spaces (the Schatten classes of order 1). Since $\T(\H)^*=\B(\H)$, the below results can be reformulated in terms of CP linear maps between algebras of all bounded operators on  separable Hilbert spaces. Nevertheless, the use of the "predual picture" is more natural for representation of our  results. The theory of CP linear maps between Banach spaces of trace class
operators  has important applications in mathematical physics, in particular, in the theory of open quantum systems, where
CP trace-preserving linear maps called \emph{quantum channels} play the role of dynamical maps (in the Schrodinger picture), while CP trace-non-increasing  linear maps
called \emph{quantum operations} are essentially used in the theory of quantum measurements \cite{H-SCI,Watrous,Wilde}.
\smallskip

For a CP linear map $\,\Phi:\T(\H_A)\rightarrow \T(\H_B)\,$ the Stinespring theorem (cf.\cite{St}) implies existence of a Hilbert space
$\mathcal{H}_E$ and an operator
$V_{\Phi}:\mathcal{H}_A\rightarrow\mathcal{H}_B\otimes\mathcal{H}_E$ such
that
\begin{equation}\label{St-rep}
\Phi(\rho)=\mathrm{Tr}_{E}V_{\Phi}\rho V_{\Phi}^{*},\quad
\rho\in\mathfrak{T}(\mathcal{H}_A),
\end{equation}
where $\Tr_E$ denotes the partial trace over $\H_E$. If $\Phi$ is trace-preserving (correspondingly, trace-non-increasing) then $V_{\Phi}$ is an isometry
(correspondingly, contraction) \cite[Ch.6]{H-SCI}.

The dual CP linear map $\,\Phi^*:\B(\H_B)\rightarrow \B(\H_A)\,$ has the corresponding representation
\begin{equation}\label{St-rep+}
\Phi^*(B)=V^*_{\Phi}[B\otimes I_E] V_{\Phi},\quad
B\in\mathfrak{B}(\mathcal{H}_B).
\end{equation}

The norm of complete boundedness ($cb\shs$-norm  in what follows) of a linear map between the algebras $\B(\H_B)$ and $\B(\H_A)$ (cf. \cite{Paul}) induces (by duality) the
norm
\begin{equation}\label{d-n-def}
\|\Phi\|_{\rm cb}\doteq\sup_{\rho\in\T(\H_{AR}),\|\rho\|_1\leq 1}\|\Phi\otimes \id_R(\rho)\|_1
\end{equation}
on the set of all linear maps between Banach spaces  $\T(\H_A)$ and $\T(\H_B)$, where $\H_R$ is a separable Hilbert space and $\H_{AR}=\H_{A}\otimes \H_{R}$. If $\Phi$ is a Hermitian preserving map
then the supremum in (\ref{d-n-def}) can be taken over the set $\S(\H_{AR})$ \cite[Ch.3]{Watrous}.


The Kretschmann-Schlingemann-Werner theorem  (the KSW-theorem in what follows) obtained in \cite{Kr&W} states that
\begin{equation*}
 \frac{\|\Phi-\Psi\|_{\rm cb}}{\sqrt{\|\Phi\|_{\rm cb}}+\sqrt{\|\Psi\|_{\rm cb}}}\leq\inf_{V_{\Phi},V_{\Psi}}\|V_{\Phi}-V_{\Psi}\|\leq\sqrt{\|\Phi-\Psi\|_{\rm cb}},
\end{equation*}
where the infimum is over all common Stinespring representations
\begin{equation}\label{c-S-r}
\Phi(\rho)=\Tr_E V_{\Phi}\rho V^*_{\Phi}\quad\textrm{and}\quad\Psi(\rho)=\Tr_E V_{\Psi}\rho V^*_{\Psi}.
\end{equation}
In the proof of the KSW theorem it is shown that the quantity $\inf_{V_{\Phi},V_{\Psi}}\|V_{\Phi}-V_{\Psi}\|$ coincides with the\emph{ Bures distance} between the maps $\Phi$ and $\Psi$ defined by the expression
\begin{equation}\label{b-dist+}
\beta(\Phi,\Psi)=\sup_{\rho\in\S(\H_{AR})} \beta\!\left(\Phi\otimes \id_R(\rho),\Psi\otimes \id_R(\rho)\right),
\end{equation}
in which $\H_R$ is a separable Hilbert space and $\beta(\cdot,\cdot)$ in the r.h.s. is the Bures distance between operators in $\T_+(\H_{BR})$ defined as
\begin{equation}\label{B-d-s}
  \beta(\rho,\sigma)=\sqrt{\|\rho\|_1+\|\sigma\|_1-2\sqrt{F(\rho,\sigma)}},
\end{equation}
where
\begin{equation}\label{fidelity}
  F(\rho,\sigma)=\|\sqrt{\rho}\sqrt{\sigma}\|^2_1
\end{equation}
is the fidelity of the operators $\rho$ and $\sigma$ \cite{H-SCI,Watrous,Wilde}. The Bures distance between CP linear maps $\Phi$ and $\Psi$
is connected to the \emph{operational fidelity} of these maps introduced in \cite{B&Co}.

The KSW theorem shows continuity of the map $V_{\Phi}\mapsto\Phi$ and selective continuity of the multi-valued
map $\Phi\mapsto V_{\Phi}$ with respect to the  $cb\shs$-norm topology on the set $\F(A,B)$ of all CP linear maps $\Phi$  from $\T(\H_A)$ to $\T(\H_B)$ and the operator norm topology on the set of Stinespring operators $V_{\Phi}$.

The $cb\shs$-norm topology is widely used in the quantum  theory, by it is too strong for description of
physical perturbations of infinite-dimensional quantum channels \cite{SCT,W-EBN}. Our aim is to obtain a version of the KSW theorem which would show continuity of the map $V_{\Phi}\mapsto\Phi$ and  selective continuity of the multi-valued
map $\Phi\mapsto V_{\Phi}$ with respect to\emph{ weaker topologies}  on the sets of CP linear maps $\Phi$ and Stinespring operators $V_{\Phi}$. A natural way to do this
is to use the operator \emph{E}-norms induced by some positive operator ${G}$ on $\H_A$ (naturally generalized
to operators between different separable Hilbert spaces, see Remark \ref{en-gen}) and the energy-constrained  $cb\shs$-norms
\begin{equation}\label{E-sn}
  \|\Phi\|^G_{\mathrm{cb},E}\doteq\sup_{\rho\in\S(\H_{AR}):\Tr {G}\rho_A\leq E}\|\Phi\otimes \id_R(\rho)\|_1,\quad E>0,\quad (\textrm{where}\; \rho_A\doteq\Tr_R\shs\rho)
\end{equation}
on the set of Hermitian-preserving
linear maps from $\T(\H_{A})$ to $\T(\H_{B})$ introduced independently in \cite{SCT} and \cite{W-EBN} (the positive operator ${G}$ is treated therein as a Hamiltonian of a quantum system $A$).\footnote{Slightly different energy-constrained $cb\shs$-norm is used in \cite{Pir}.} If ${G}$ is a discrete unbounded operator (see Def.\ref{D-H}) then the topology
generated by any of the norms (\ref{E-sn}) on bounded subsets of $\F(A,B)$ coincides with the strong  convergence topology generated  by the family of seminorms $\Phi\mapsto\|\Phi(\rho)\|_1$, $\rho\in\T(\H_A)$ \cite[Proposition 3]{SCT}.\footnote{This topology is a restriction
to the set $\F(A,B)$ of the strong operator topology on the set of all linear maps from $\T(\H_A)$ to $\T(\H_B)$. The strong convergence of a sequence $\{\Phi_n\}\subset\F(A,B)$  to a  map $\Phi_0$  means that
$\lim_{n\rightarrow\infty}\Phi_n(\rho)=\Phi_0(\rho)\,\textup{ for all }\rho\in\T(\H_A)$.}

Following \cite{CID} introduce the \emph{energy-constrained Bures distance}
\begin{equation}\label{ec-b-dist}
\beta_E^G(\Phi,\Psi)=\sup_{\rho\in\S(\H_{AR}):\Tr {G}\rho_A\leq E} \beta(\Phi\otimes \id_R(\rho),\Psi\otimes \id_R(\rho)), \quad E>0,
\end{equation}
between CP linear maps $\Phi$  and $\Psi$ from $\T(\H_A)$ to $\T(\H_B)$, where $\beta(\cdot,\cdot)$ in the r.h.s. is the Bures distance between operators in $\T_+(\H_{BR})$ defined in (\ref{B-d-s}) and $\H_R$ is a separable Hilbert space.
\begin{remark}\label{ec-b-dist-r} The infimum in (\ref{ec-b-dist}) can be taken only over pure states $\rho\in\S(\H_{AR})$.
This follows from the freedom of choice of $R$, which implies possibility to purify any mixed state in $\S(\H_{AR})$ by extending system $R$. We have only to note that
the Bures distance between operators in $\T_{+}(\H_{XY})$ defined in (\ref{B-d-s}) does not increase under partial trace: $\beta(\rho,\sigma)\geq \beta(\rho_X,\sigma_X)$ for any $\rho$ and $\sigma$ in $\T_{+}(\H_{XY})$ \cite{H-SCI,Watrous,Wilde}. \smallskip
\end{remark}
The distance $\beta_E^G(\Phi,\Psi)$ turns out to be extremely useful in quantitative continuity analysis of capacities of energy-constrained  infinite-dimensional quantum channels \cite[Theorem 2]{CID}.  By using the well known relations between the trace norm and the Bures distance (\ref{B-d-s}) one can show that for any $E>0$ the distance $\beta_E^G(\Phi,\Psi)$ generates the same topology on bounded subsets of $\F(A,B)$ as any of the
energy-constrained  $cb\shs$-norms (\ref{E-sn}). The results of calculation of $\beta_E^G(\Phi,\Psi)$
for real quantum channels can be found in \cite{Nair}.\smallskip

Now we can formulate the $E$-version of  KSW-theorem.\smallskip

\begin{theorem}\label{KSW-E} \emph{Let ${G}$ be a positive semidefinite densely defined operator on $\H_A$ satisfying condition (\ref{H-cond}) and $E>0$. Let $\|\cdot\|^G_{\mathrm{cb},E}$ and $\|\cdot\|_E^G$
be, respectively,  the energy-constrained  $cb\shs$-norm and the operator $E\textrm{-}$norm induced  by ${G}$. For any CP linear maps $\,\Phi$ and $\,\Psi$ from $\,\T(\H_A)$ to $\,\T(\H_B)$ the following inequalities hold
\begin{equation}\label{KSW-rel}
 \frac{\|\Phi-\Psi\|^G_{\mathrm{cb},E}}{\sqrt{\|\Phi\|^G_{\mathrm{cb},E}}+\sqrt{\|\Psi\|^G_{\mathrm{cb},E}}}\leq\inf_{V_{\Phi},V_{\Psi}}
 \|V_{\Phi}-V_{\Psi}\|_E^G\leq\sqrt{\|\Phi-\Psi\|^G_{\mathrm{cb},E}},
\end{equation}
where the infimum is over all common Stinespring representation (\ref{c-S-r}). The quantity $\,\inf_{V_{\Phi},V_{\Psi}}\|V_{\Phi}-V_{\Psi}\|_E^G$
coincides with the energy-constrained Bures distance $\beta_E^G(\Phi,\Psi)$ defined in (\ref{ec-b-dist}). The infimum in (\ref{KSW-rel}) is attainable.}
\end{theorem}\medskip

\emph{Proof.} We will follow the proof of the KSW theorem (given in \cite{Kr&W})
with necessary modifications concerning the use of the energy-constrained $cb\shs$-norms and
the operator $E$-norms (instead of the ordinary  $cb\shs$-norm and the operator norm).

To prove the first inequality in (\ref{KSW-rel}) assume that $\rho$ is a state in $\S(\H_{AR})$ such that $\Tr {G}\rho_A\leq E$. For a given common Stinespring representation (\ref{c-S-r}) we have
$$
\begin{array}{c}
\|(\Phi-\Psi)\otimes \id_R(\rho)\|_1\leq \|V_{\Phi}\otimes I_R \cdot\rho\cdot V^*_{\Phi}\otimes I_R-V_{\Psi}\otimes I_R \cdot\rho\cdot V^*_{\Psi}\otimes I_R\|_1\\\\
\leq\|(V_{\Phi}-V_{\Psi})\otimes I_R \cdot\rho\cdot V^*_{\Phi}\otimes I_R\|_1+\|V_{\Psi}\otimes I_R \cdot\rho\cdot(V^*_{\Phi}-V^*_{\Psi})\otimes I_R\|_1\\\\
\leq\|(V_{\Phi}-V_{\Psi})\otimes I_R\|_E^{G\otimes I_R} \|V_{\Phi}\otimes I_R\|_E^{G\otimes I_R}+\|(V_{\Phi}-V_{\Psi})\otimes I_R\|_E^{G\otimes I_R}\|V_{\Psi}\otimes I_R\|_E^{G\otimes I_R}
\\\\
\leq\|V_{\Phi}-V_{\Psi}\|_E^G \|V_{\Phi}\|_E^G+\|V_{\Phi}-V_{\Psi}\|_E^G \|V_{\Psi}\|_E^G.
\end{array}
$$
The first and the second inequalities follow from the properties of the trace norm (the non-increasing under partial trace and the triangle inequality), the third inequality follows from Proposition \ref{bp-en-0}D, the last one -- from Proposition \ref{pen-tp}A. By noting that $[\|V_{\Phi}\|^G_E]^2=\|\Phi\|^G_{\mathrm{cb},E}$ and
$[\|V_{\Psi}\|^G_E]^2=\|\Psi\|^G_{\mathrm{cb},E}$ we obtain the first inequality in (\ref{KSW-rel}).

To prove the second inequality in (\ref{KSW-rel}) note that
$\beta_E^G(\Phi,\Psi)\leq\sqrt{\|\Phi-\Psi\|^G_{\mathrm{cb},E}}$.
This follows from the inequality $\beta(\rho,\sigma)\leq\sqrt{\|\rho-\sigma\|_1}$ valid for any $\rho$ and $\sigma$ in $\T_+(\H)$, which is easily proved by using the inequality $\Tr(\sqrt{\rho}-\sqrt{\sigma})^2\leq \|\rho-\sigma\|_1$ (see the proof of Lemma 9.2.3 in \cite{H-SCI}).
So, it suffices to show that
\begin{equation}\label{r-r}
\inf_{V_{\Phi},V_{\Psi}}\|V_{\Phi}-V_{\Psi}\|_E^G=\beta_E^G(\Phi,\Psi).
\end{equation}

Denote by $\alpha^G_E(\Phi,\Psi)$ the l.h.s. of (\ref{r-r}).
Let $\C^s_{{G},E}$ be the subset of $\S(\H_A)$ determined by the inequality $\,\Tr {G}\rho\leq E\,$ and $\N(\Phi,\Psi)=\bigcup V_{\Phi}^*V_{\Psi}$, where the union is over \emph{all} common Stinespring representations (\ref{c-S-r}). Then by using definition (\ref{ec-on}) we obtain
\begin{equation}\label{beta-e}
\alpha^G_E(\Phi,\Psi)=\inf_{N\in\N(\Phi,\Psi)} \sup_{\rho\in\C^s_{{G},E}}\sqrt{\Tr\shs\Phi(\rho)+\Tr\Psi(\rho)-2\Re\shs\Tr N\rho}.
\end{equation}
Following the proof of Theorem 1 in \cite{Kr&W} show that $\N(\Phi,\Psi)$ coincides with the set
$$
\M(\Phi,\Psi)\doteq \left\{V_{\Phi} ^*(I_B\otimes C)V_{\Psi}\,|\,C\in \B(\H_E), \|C\|\leq 1\right\},
$$ defined via some \emph{fixed} common Stinespring representation (\ref{c-S-r}). It will imply, in particular, that
$\M(\Phi,\Psi)$ does not depend on this representation.

To show that $\M(\Phi,\Psi)\subseteq\N(\Phi,\Psi)$ it suffices to find for any contraction $C\in\B(\H_E)$ a common Stinespring representation for  $\Phi$ and $\Psi$ with the operators
$\tilde{V}_{\Phi}$ and $\tilde{V}_{\Psi}$ from $\H_A$ to $\H_B\otimes\H_{\tilde{E}}$ such that $\tilde{V}^*_{\Phi}\tilde{V}_{\Psi}=V_{\Phi} ^*(I_B\otimes C)V_{\Psi}$.

Let $\H_{\tilde{E}}=\H^1_E\oplus\H^2_E$, where $\H^1_E$ and $\H^2_E$ are copies of $\H_E$. For given $C$ define the operators
$\tilde{V}_{\Phi}$ and $\tilde{V}^C_{\Psi}$ from $\H_A$ into $\H_B\otimes(\H_{E_1}\oplus\H_{E_2})=\H_B\otimes\H_{E_1}\oplus\H_B\otimes\H_{E_2}$ by setting
\begin{equation}\label{bar-v}
\tilde{V}_{\Phi}|\varphi\rangle=V_{\Phi}|\varphi\rangle\oplus|0\rangle,\quad
\tilde{V}^C_{\Psi}|\varphi\rangle=(I_B\otimes C)V_{\Psi}|\varphi\rangle\oplus \left(I_B\otimes\sqrt{I_{E}-C^*C}\right)V_{\Psi}|\varphi\rangle
\end{equation}
for any $\varphi\in\H_A$,
where we assume that the operators $V_{\Phi}$ and $V_{\Psi}$ act from $\H_A$ to $\H_B\otimes \H^1_E$ and
$\H_B\otimes \H^2_E$ correspondingly, while the contraction $C$ acts from $\H^2_E$ to $\H^1_E$.  It is easy to see that
the operators $\tilde{V}_{\Phi}$ and $\tilde{V}^C_{\Psi}$ form a common Stinespring representation for the maps $\Phi$ and $\Psi$ with the required property.

To prove that $\N(\Phi,\Psi)\subseteq\M(\Phi,\Psi)$ take any common Stinespring representation for the maps $\Phi$ and $\Psi$ with the operators
$\tilde{V}_{\Phi}$ and $\tilde{V}_{\Psi}$ from $\H_A$ to $\H_B\otimes\H_{\tilde{E}}$. By Theorem 6.2.2 in \cite{H-SCI} there exist partial isometries $W_{\Phi}$ and $W_{\Psi}$ from $\H_E$ to $\H_{\tilde{E}}$ such that $\tilde{V}_{\Phi}=(I_B\otimes W_{\Phi})V_{\Phi}$ and $\tilde{V}_{\Psi}=(I_B\otimes W_{\Psi})V_{\Psi}$. So, $\tilde{V}^*_{\Phi}\tilde{V}_{\Psi}=V^*_{\Phi}(I_B\otimes W^*_{\Phi}W_{\Psi})V_{\Psi}\in \M(\Phi,\Psi)$, since $\|W^*_{\Phi}W_{\Psi}\|\leq 1$.

Since $\N(\Phi,\Psi)=\M(\Phi,\Psi)$, the infimum in (\ref{beta-e}) can be taken over the set $\M(\Phi,\Psi)$. This implies
\begin{equation}\label{beta-d}
\begin{array}{rl}
\alpha^G_E(\Phi,\Psi)&
=\displaystyle\inf_{C\in\B_1(\H_E)}\sup_{\rho\in\C^s_{{G},E}}\sqrt{\Tr\shs\Phi(\rho)+\Tr\Psi(\rho)-2\Re\shs\Tr V_{\Phi}^*(I_B\otimes C)V_{\Psi}\rho}\\
&=\displaystyle\sup_{\rho\in\C^s_{{G},E}}\inf_{C\in\B_1(\H_E)} \sqrt{\Tr\shs\Phi(\rho)+\Tr\Psi(\rho)-2\Re\shs\Tr V_{\Phi} ^*(I_B\otimes C)V_{\Psi}\rho}\\
&=\displaystyle\sup_{\rho\in\C^s_{{G},E}}\sqrt{\shs\Tr\shs\Phi(\rho)+\Tr\Psi(\rho)-2\sup_{C\in\B_1(\H_E)}|\Tr V_{\Phi} ^*(I_B\otimes C)V_{\Psi}\rho|},
\end{array}
\end{equation}
where the possibility to change the
order of the optimization follows from Ky Fan's minimax theorem \cite{Simons} and the $\sigma$-weak compactness of the unit ball $\B_1(\H_E)$ of $\B(\H_E)$ \cite{B&R}. It is easy to see that
\begin{equation}\label{a-eq}
\sup_{C\in\B_1(\H_E)}|\Tr V_{\Phi} ^*(I_B\otimes C)V_{\Psi}\rho|=\!\!\sup_{C\in\B_1(\H_E)}|\langle V_{\Phi}\otimes I_R\shs\varphi|I_{BR}\otimes C |V_{\Psi}\otimes I_R \shs\varphi\rangle|,\!
\end{equation}
where $\varphi$ is a purification of $\rho$, i.e. a vector in $\H_{A}\otimes\H_R$ such that $\Tr_R|\varphi\rangle\langle\varphi|=\rho$.

Since for \emph{any} common Stinespring representation (\ref{c-S-r}) and any purification $\varphi$ of a state $\rho$ the vectors $V_{\Phi}\otimes I_R\shs|\varphi\rangle$ and $V_{\Psi}\otimes I_R\shs|\varphi\rangle$ in $\H_{BER}$ are purifications of the operators
$\Phi\otimes \id_R(|\varphi\rangle\langle\varphi|)$ and $\Psi\otimes \id_R(|\varphi\rangle\langle\varphi|)$ in $\T(\H_{BR})$, by using the relation $\N(\Phi,\Psi)=\M(\Phi,\Psi)$ proved before and Uhlmann's theorem \cite{Uhlmann,Wilde} it is easy to show that the square of the r.h.s. of (\ref{a-eq}) coincides with the fidelity of these operators defined in (\ref{fidelity}).
Note also that $\Tr\shs\Phi\otimes \id_R(\sigma)=\Tr\shs\Phi(\sigma_A)$ and  $\Tr\Psi\otimes \id_R(\sigma)=\Tr\Psi(\sigma_A)$ for any state $\sigma$ in $\S(\H_{AR})$. By Remark \ref{ec-b-dist-r} these observations and (\ref{beta-d}) imply that $\alpha^G_E(\Phi,\Psi)=\beta_E^G(\Phi,\Psi)$, i.e. that (\ref{r-r}) holds.

The last assertion  can be derived from the attainability of the infimum in the first line in (\ref{beta-d}) which follows from the $\sigma$-weak compactness of the unit ball $\B_1(\H_E)$. $\square$\medskip

Theorem \ref{KSW-E} shows continuity of the map $V_{\Phi}\mapsto\Phi$ and  selective continuity of the multi-valued
map $\Phi\mapsto V_{\Phi}$ with respect to the energy-constrained $cb\shs\textrm{-}$norm on the set of CP linear maps $\Phi$ and the operator \emph{E}-norm on the set of Stinespring operators $V_{\Phi}$. Its basic assertion is the equality
\begin{equation}\label{ksw-rel-+}
\beta_E^G(\Phi,\Psi)=\inf_{V_{\Phi},V_{\Psi}} \|V_{\Phi}-V_{\Psi}\|_E^G.
\end{equation}
Some difficulty of applying Theorem \ref{KSW-E} is related to the fact that the infimum in
(\ref{ksw-rel-+}) is over all common Stinespring representation (\ref{c-S-r}). But by using the constructions from the proof of this theorem one can
obtain its versions which are more convenient for applications, in particular, for analysis of converging sequences of CP linear maps. \smallskip

\begin{theorem}\label{KSW-E+} \emph{Let ${G}$ be a positive semidefinite densely defined operator on $\H_A$ satisfying condition (\ref{H-cond}), $\beta_E^G$ and $\|\!\cdot\!\|_E^G$
be, respectively,  the energy-constrained  Bures distance and the operator $E\textrm{-}$norm induced by ${G}$. Let $\,\Phi$ be a CP linear map from $\T(\H_A)$ to $\T(\H_B)$.}\smallskip

\noindent A) \emph{There is  a Stinespring representation of $\,\Phi$ with the operator $\,V'_{\Phi}:\H_A\rightarrow\H_B\otimes\H_{E'}$ such that
\begin{equation}\label{ksw-rel+}
\beta_E^G(\Phi,\Psi)=\inf_{V_{\Psi}}
 \|V'_{\Phi}-V_{\Psi}\|_E^G,
\end{equation}
for any CP linear map $\Psi:\T(\H_A)\rightarrow\T(\H_B)$, where the infimum is over all Stinespring representations of $\,\Psi$ with the same environment space $\,\H_{E'}$. The infimum in (\ref{ksw-rel+}) is attainable.}\smallskip

\noindent B) \emph{If $G$~is an unbounded discrete operator (Def.\ref{D-H}) and  $\,V_{\Phi}:\H_A\rightarrow\H_B\otimes\H_{E}$  is the operator from a given Stinespring representation of $\,\Phi$ such that $\,\dim\H_{E}=+\infty\,$ then
$$
\beta_E^G(\Phi,\Psi)\leq \inf_{V_{\Psi}}\|V_{\Phi}-V_{\Psi}\|_E^G\leq 2\beta_E^G(\Phi,\Psi),
$$
for any CP linear map $\,\Psi:\T(\H_A)\rightarrow\T(\H_B)$, where the infimum is over all Stinespring representations of $\,\Psi$ with the same environment space $\H_{E}$.}
\end{theorem}\medskip

\emph{Proof.} If  $\,V_{\Phi}:\H_A\rightarrow\H_B\otimes\H_{E}$ is the operator from a Stinespring representation of $\,\Phi$ such that $\,\dim\H_E=+\infty\,$ then,
since any separable Hilbert space can be isometrically embedded into $\H_E$, we may assume that any CP linear map $\Psi:\T(\H_A)\rightarrow\T(\H_B)$ has a Stinespring representation with the same environment space $\H_E$. Denote by $V_{\Psi}$ the Stinespring operator of  $\Psi$ in this representation.
Let $\tilde{V}_{\Phi}$ and $\tilde{V}^C_{\Psi}$ be the operators from $\H_A$ into $\H_B\otimes(\H^1_E\oplus\H^2_E)=(\H_B\otimes\H^1_E)\oplus(\H_B\otimes\H^2_E)$  defined by formulae (\ref{bar-v}), where $\H^1_E$ and $\H^2_E$ are copies of $\H_E$ and $C$ is a contraction in $\B(\H_E)$.  The arguments from the proof of  Theorem \ref{KSW-E} show that
$\,\beta_E^G(\Psi,\Phi)=\|\tilde{V}^{C_0}_{\Psi}-\tilde{V}_{\Phi}\|_E^G$ for some $C_0\in\B(\H_E)$ depending on $\Phi$ and $\Psi$.
So, to obtain assertion A it suffices to take $\tilde{V}_{\Phi}$ in the role of $V'_{\Phi}$.

To prove assertion B we will use the above operators $\tilde{V}_{\Phi}$ and $\tilde{V}^{C_0}_{\Psi}$  as follows. Assume first that the operator $C_0$ is nondegenerate, i.e. $\ker C_0=\{0\}$. Let $U$ be the isometry from the polar decomposition of $C_0$, i.e. $C_0=U|C_0|$. Since $\,\|\tilde{V}^{C_0}_{\Psi}-\tilde{V}_{\Phi}\|_E^G=\beta_E^G(\Psi,\Phi)$, it follows from Proposition \ref{bp-en-0}C that
\begin{equation}\label{s-ub}
\!\!\|(I_B\otimes C_0)V_{\Psi}-V_{\Phi}\|_E^G\leq\beta_E^G(\Psi,\Phi)\quad\textrm{and}\quad \left\|\left(I_B\otimes\sqrt{I_{E}-|C_0|^2}\right)\!V_{\Psi}\right\|_E^G\!\leq\beta_E^G(\Psi,\Phi)\!
\end{equation}
Hence the triangle inequality and Proposition \ref{bp-en-0}B imply that
\begin{equation}\label{s-ub+}
\begin{array}{c}
\|(I_B\otimes U)V_{\Psi}-V_{\Phi}\|_E^G\leq \|(I_B\otimes C_0)V_{\Psi}-V_{\Phi}\|_E^G\\\\+\|(I_B\otimes C_0)V_{\Psi}-(I_B\otimes U)V_{\Psi}\|_E^G
\leq \beta_E^G(\Psi,\Phi)+\|I_B\otimes (I_E-|C_0|)V_{\Psi}\|_E^G.
\end{array}
\end{equation}
Since $C_0$ is a contraction, by using Proposition \ref{bp-en-0}B and the second inequality in (\ref{s-ub}) we obtain
$$
\|I_B\otimes (I_E-|C_0|)V_{\Psi}\|_E^G\leq\|I_B\otimes (I_E-|C_0|^2)V_{\Psi}\|_E^G\leq\|I_B\otimes \sqrt{I_E-|C_0|^2}V_{\Psi}\|_E^G\leq\beta_E^G(\Psi,\Phi)
$$
Thus, it follows from (\ref{s-ub+}) that $\|(I_B\otimes U)V_{\Psi}-V_{\Phi}\|_E^G\leq2\beta_E^G(\Psi,\Phi)$. Since $U$ is an isometry,
$(I_B\otimes U)V_{\Psi}$ is a Stinespring operator for $\Psi$.\smallskip

Since $G$~is a discrete unbounded operator on~$\H_A$, the set $\C^s_{{G},E}$ of states $\rho$ in $\mathfrak{S}(\H_A)$ such that $\Tr {G}\rho\leq E$ is compact by Lemma \ref{Comp}. Using this and taking into account the continuity of the expression under the square root in the first line in~(\ref{beta-d}) as a function on the
Cartesian product of the set $\C^s_{{G},E}$ and the set  $\mathfrak{B}_1(\H_E)$ equipped with the weak operator topology, it is easy to show that the first infimum in~(\ref{beta-d}) can be taken over the dense subset of $\B_1(\H_E)$ consisting of non-degenerate operators. This allows to omit the assumption $\ker C_0=\{0\}$.
$\square$ \smallskip

If $\{V_n\}$ is a sequence of operators from $\H_A$ to $\H_B\otimes\H_E$ converging  to an operator $V_0:\H_A\rightarrow\H_B\otimes\H_E$ w.r.t. the norm $\|\!\cdot\!\|^G_{E}$ then the first inequality in (\ref{KSW-rel}) implies that the sequence of CP maps $\,\Phi_n(\rho)=\Tr_E V_n\rho V^*_n\,$ converges to the map $\,\Phi_0(\rho)=\Tr_E V_0\rho V^*_0\,$ w.r.t. the norm $\|\!\cdot\!\|^G_{\mathrm{cb},E}$  and for each $n$ the following inequalities hold
$$
\|\Phi_n-\Phi_0\|^G_{\mathrm{cb},E}\leq\beta^G_{E}(\Phi_n,\Phi_0)\!\left[\sqrt{\|\Phi_n\|^G_{\mathrm{cb},E}}+\sqrt{\|\Phi_0\|^G_{\mathrm{cb},E}}\right]
\leq\|V_{n}-V_{0}\|_E^G\!\left[\|V_{n}\|_E^G+\|V_{0}\|_E^G\right].
$$

Theorem \ref{KSW-E+} allows to describe all  sequences of CP linear maps converging w.r.t. the energy-constrained $cb\shs\textrm{-}$norm. \smallskip

\begin{corollary}\label{KSW-E+c} \emph{Let $\,\{\Phi_n\}$ be a sequence of CP linear maps from $\,\T(\H_A)$ to $\T(\H_B)$ converging to a CP linear map $\,\Phi_0$ with respect to the norm $\|\!\cdot\!\|^G_{\mathrm{cb},E}$.}

\smallskip

\noindent A) \emph{There exist a separable Hilbert space $\H_{E'}$ and a sequence $\{V_n\}$ of operators from $\H_A$ into $\H_B\otimes\H_{E'}$ converging to an operator $V_0$ with respect to the norm $\|\!\cdot\!\|^G_{E}$ such that $\,\Phi_n(\rho)=\Tr_{E'} V_n\rho V^*_n\,$ for all $\,n\geq0$ and}
$$
\|V_{n}-V_{0}\|_E^G=\beta_{E}^G(\Phi_n,\Phi_0)\leq \sqrt{\|\Phi_n-\Phi_0\|^G_{\mathrm{cb},E}}.
$$

\noindent B) \emph{If $G$~is an unbounded discrete operator (Def.\ref{D-H}) and $\,V_0:\H_A\rightarrow\H_B\otimes\H_E$ is the operator from a given Stinespring representation of the map $\,\Phi_0$ such that $\,\dim \H_E=+\infty$,  then for any  $\varepsilon>0$ there exists  a sequence $\{V_n\}$ of operators from $\H_A$ into $\H_B\otimes\H_E$ converging to the operator $V_0$ with respect to the norm $\|\!\cdot\!\|_{E}^G$ such that $\,\Phi_n(\rho)=\Tr_E V_n\rho V^*_n\,$ for all $\,n>0$ and}
\begin{equation}\label{2-est}
\|V_{n}-V_{0}\|_E^G\leq2\beta^G_{E}(\Phi_n,\Phi_0)+\varepsilon\leq 2\sqrt{\|\Phi_n-\Phi_0\|^G_{\mathrm{cb},E}}+\varepsilon.
\end{equation}
\end{corollary}
\smallskip

Factor $"2"$ in (\ref{2-est}) is a cost of the possibility to take the sequence $\{V_n\}$ of Stinespring operators representing the sequence $\{\Phi_n\}$ for \emph{given} $\H_E$ and $V_0:\H_A\rightarrow\H_B\otimes\H_E$.\smallskip

If the  operator ${G}$ is discrete and  unbounded (Def.\ref{D-H}) then the norm $\|\!\cdot\!\|^G_{\mathrm{cb},E}$\break
generates the strong convergence topology on bounded subsets of the set $\F(A,B)$ of all CP linear maps from $\T(\H_A)$ to $\T(\H_B)$ (by Proposition 3 in \cite{SCT}), while the norm $\|\!\cdot\!\|^G_{E}$ generates the strong operator topology on subsets of linear maps from $\H_A$ to $\H_B\otimes\H_E$ bounded by the operator norm (by Proposition \ref{ec-on-p2}B). Thus, in this case Corollary \ref{KSW-E+c} gives  representation of bounded strongly converging
sequences of CP linear maps via  strongly converging sequence of Stinespring operators. For sequences of quantum channels such representation is obtained in \cite{CSR}.

\section{Operator \emph{E}-norms for unbounded operators}

In this section we will extend the operator \emph{E}-norms  to unbounded operators.  We will assume that ${G}$ is
a positive semidefinite \emph{unbounded}\footnote{If ${G}$ is a bounded operator then the norm $\|\!\cdot\!\|^{G}_E$ is equivalent to the operator norm by Proposition 2A.} operator on $\H$ with dense domain satisfying condition (\ref{H-cond}). The case of discrete type operator ${G}$ will be considered separately after formulations of general results.

Speaking about extension of the operator \emph{E}-norms to unbounded operators we may restrict attention to
linear operators on $\H$ relatively bounded  w.r.t. the operator $\sqrt{G}$, i.e. linear operators $A$ defined on $\mathcal{D}(\sqrt{G})$ such that
\begin{equation}\label{rb-rel}
\|A\varphi\|^2\leq a^2\|\varphi\|^2+b^2\|\sqrt{G}\varphi\|^2,\quad \forall \varphi\in\mathcal{D}(\sqrt{G}),
\end{equation}
for some  nonnegative numbers $a$ and $b$ (depending on $A$ but not depending on $\varphi$)\cite{Kato}. Such operators are briefly called $\sqrt{G}$-\emph{bounded}.
Indeed, it is easy to see that the r.h.s. of (\ref{ec-on-b}) is finite for any $\sqrt{G}$-bounded operator $A$ and all $E>0$.
The following lemma contains the converse statement (in strengthened form).

\begin{lemma}\label{new-l} \emph{Let $A$ be a linear operator on $\H$ such that $\D(A)\supseteq\D(\sqrt{G})$. If
the quantity $\|A\|_E^G$ defined in (\ref{ec-on-b}) is finite for some $E>0$ then
\begin{itemize}
  \item the function $E\mapsto \left[\|A\|^{G}_E\right]^2$ is finite and concave on $\,\mathbb{R}_{+}$;
  \item the operator $A$ is $\sqrt{G}$-bounded.
\end{itemize}}
\end{lemma}

\emph{Proof.} Consider the set
$\S^{\rm f}_{\!G}\doteq \{ \rho\in \S(\H)\,|\,\Tr {G}\rho<+\infty, \rank \rho<+\infty\}$.
A state $\rho$ belongs to this set if and only if it has a
finite decomposition
\begin{equation}\label{p-s-d}
  \rho=\sum_i |\varphi_i\rangle\langle\varphi_i|,
\end{equation}
where  $\{\varphi_i\}$ is a set of  vectors in $\D(\sqrt{G})$. Moreover, any
such decomposition of $\rho$ consists of  vectors in $\D(\sqrt{G})$.

For any state $\rho$ in $\S^{\rm f}_{\!G}$ with representation (\ref{p-s-d}) define the operator $A\rho A^*$ as follows
\begin{equation}\label{ab-d-p}
  A\rho A^*\doteq\sum_i|\alpha_i\rangle\langle\alpha_i|,\quad\textrm{where}\quad |\alpha_i\rangle=A|\varphi_i\rangle.
\end{equation}
By using Schrodinger's mixture theorem (see \cite[Ch.8]{B&Z}) it is easy to show that the r.h.s. of (\ref{ab-d-p}) does not depend on  representation (\ref{p-s-d}). This implies that
\begin{equation}\label{a-fun}
\rho\mapsto\Tr A\rho A^*\doteq\sum_i\|A\varphi_i\|^2
\end{equation}
is an affine function on $\S^{\rm f}_{\!G}$. \smallskip

Thus, to prove the first assertion of the lemma it suffices to show that
\begin{equation}\label{ec-on-f}
 \|A\|^{G}_E=\sup_{\rho\in\S^{\rm f}_{\!G}:\Tr {G}\rho\leq E}\sqrt{\Tr A\rho A^*}
\end{equation}
for any $E>0$. This equality means that the supremum in the r.h.s. of (\ref{ec-on-f}) can be taken only over
pure states $\rho$ in $\S^{\rm f}_{\!G}$ such that $\Tr {G}\rho\leq E$.
Since the function (\ref{a-fun}) is affine, this can be shown easily by using the fact that all the extreme points of the convex set of states $\rho$ such that $\Tr {G}\rho\leq E$ are pure states \cite{W&Sh}.

The second assertion of the lemma is derived from the first one, since the concavity of the function $\,E\mapsto \left[\|A\|^{G}_E\right]^2$  implies existence of  numbers $a$ and $b$ such that
$$
[\|A\|^{G}_E]^2\leq a^2+b^2 E\quad \forall E>0.\;\square
$$
\smallskip

Thus, in what follows we will consider the operator $E$-norm $\|\cdot\|_E^G$ defined by formula (\ref{ec-on-b}) on the set of all $\sqrt{G}$-bounded operators. Below we will show that the quantity $\Tr A\rho A^*$ can be defined correctly (without using the notion of adjoint operator)
for any $\sqrt{G}$-bounded operator $A$ and any state $\rho$ with finite $\Tr G\rho$ (not only for a finite rank state as in the proof of Lemma \ref{new-l}). This will allows to show that the operator $E$-norm $\|A\|_E^G$ of any $\sqrt{G}$-bounded operator $A$
can be also defined by formula (\ref{ec-on}).\smallskip

Denote by $\Pi_{\sqrt{G}}(A)$ the set of all pairs $(a,b)$ for which
(\ref{rb-rel}) holds. It is easy to see that $\Pi_{\sqrt{G}}(A)$ is a closed subset of $\mathbb{\mathbb{R}}^2_+$. The  $\sqrt{G}$-bound of $A$ (denoted by $b_{\sqrt{G}}(A)$ in what follows) is defined as
$$
b_{\sqrt{G}}(A)=\inf\left\{b\,|\,(a,b)\in\Pi_{\sqrt{G}}(A)\right\}.
$$
If $\,b_{\sqrt{G}}(A)=0\,$  then $A$ is called $\sqrt{G}$-infinitesimal operator (infinitesimally bounded w.r.t. $\sqrt{G}$). These notions are widely used in
the modern operator theory, in particular, in analysis of perturbations of unbounded operators on a Hilbert space \cite{Kato,R&S+,BS}.
\smallskip

We will use the following simple lemmas.\footnote{I would be grateful for direct references to these results.}\smallskip

\begin{lemma}\label{p-l-1}
\emph{If $A$ is a $\sqrt{G}$-bounded operator on $\H$ then for any separable Hilbert space $\K$  the operator $A\otimes I_{\K}$ naturally defined on the set $\,\D(\sqrt{G})\otimes\K$  has a unique linear $\sqrt{G}\otimes I_{\K}$-bounded extension to the set $\,\D(\sqrt{G}\otimes I_{\K})$.\footnote{$\D(\sqrt{G})\otimes \K$ is the linear span of all the vectors $\varphi\otimes\psi$, where $\varphi\in\D(\sqrt{G})$ and $\psi\in\K$.} This extension (also denoted by $A\otimes I_{\K}$) has the following property
\begin{equation}\label{s-prop}
 A\otimes I_{\K}\!\left(\sum_{i}|\varphi_i\rangle\otimes|\psi_i\rangle\right)=\sum_{i}A|\varphi_i\rangle\otimes|\psi_i\rangle
\end{equation}
for any countable sets $\{\varphi_i\}\subset\D(\sqrt{G})$ and $\{\psi_i\}\subset\K$ such that $\sum_{i}\|\sqrt{G}\varphi_i\|^2<+\infty$, $\sum_{i}\|\varphi_i\|^2<+\infty$ and $\langle\psi_i|\psi_j\rangle=\delta_{ij}$, which implies that $\,\Pi_{\sqrt{G}\otimes I_{\K}}(A\otimes I_{\K})=\Pi_{\sqrt{G}}(A)$.}\smallskip
\end{lemma}

\emph{Proof.} For any $E>0$ the linear spaces $\D(\sqrt{G})$ and $\D(\sqrt{G}\otimes I_{\K})$ equipped, respectively, with the inner products
$$
\langle\varphi|\psi\rangle_E^G=\langle\varphi|\psi\rangle+\langle\varphi|G|\psi\rangle/E\quad\textrm{and}\quad  \langle\eta|\theta\rangle_E^{G\otimes I_{\K}}=\langle\eta|\theta\rangle+\langle\eta|G\otimes I_{\K}|\theta\rangle/E
$$
are  Hilbert spaces \cite{R&S}. Denote the first space by $\H^G_E$. Then it is easy to see that the second
space coincides with the Hilbert space $\H^G_E\otimes\K$.  Since the operator $A$ is bounded as an operator from $\H^G_E$ into $\H$ the operator
$A\otimes I_{\K}$ defined on $\,\D(\sqrt{G})\otimes\K$  is uniquely extended to a bounded operator from $\H^G_E\otimes\K$ into $\H\otimes\K$. Since the linear spaces $\H^G_E\otimes\K$ and $\D(\sqrt{G}\otimes I_{\K})$ coincide, this extension is a $\sqrt{G}\otimes I_{\K}$-bounded linear operator on $\H\otimes\K$.\smallskip

Property (\ref{s-prop}) follows from continuity of the operator $A\otimes I_{\K}:\H^G_E\otimes\K\rightarrow \H\otimes\K$. \smallskip

Any vector $\eta$ in $\D(\sqrt{G}\otimes I_{\K})$ can be  represented as
$|\eta\rangle=\sum_{i}|\varphi_i\rangle\otimes|\psi_i\rangle$,
where $\{\varphi_i\}\subset\D(\sqrt{G})$ and $\{\psi_i\}\subset\K$ are collections of vectors such that $\sum_{i}\|\sqrt{G}\varphi_i\|^2<+\infty$, $\sum_{i}\|\varphi_i\|^2<+\infty$ and $\langle\psi_i|\psi_j\rangle=\delta_{ij}$. By using property (\ref{s-prop}) we obtain
$$
\|A\otimes I_{\K}\eta\|^2=\sum_{i}\|A\varphi_i\|^2\leq a^2\sum_{i}\|\varphi_i\|^2+b^2\sum_{i}\|\sqrt{G}\varphi_i\|^2=a^2\|\eta\|^2+b^2\|\sqrt{G}\otimes I_{\K}\eta\|^2
$$
for any $(a,b)\in\Pi_{\sqrt{G}}(A)$. This implies $\,\Pi_{\sqrt{G}}(A)\subseteq\Pi_{\sqrt{G}\otimes I_{\K}}(A\otimes I_{\K})$, and hence
$\,\Pi_{\sqrt{G}}(A)=\Pi_{\sqrt{G}\otimes I_{\K}}(A\otimes I_{\K})$, since the converse inclusion is obvious. $\square$\medskip

\begin{remark}\label{w-c} Property (\ref{s-prop}) implies that
$$
(A\otimes I_{\K})(I_{\H}\otimes W)|\varphi\rangle=(I_{\H}\otimes W)(A\otimes I_{\K})|\varphi\rangle
$$
for any $\varphi\in\D(\sqrt{G}\otimes I_{\K})$ and a  partial isometry $W\in\B(\K)$ s.t.
$I_{\H}\otimes W^{*}W |\varphi\rangle=|\varphi\rangle$.
\end{remark}\medskip

\begin{lemma}\label{p-l-2} \emph{For any $\sqrt{G}$-bounded operators $A$ and $B$ on $\H$
the affine function $\shs\rho\mapsto A\rho B^*\in\T(\H)$ is well defined on the set $\,\T^+_{G}\doteq \{\rho\in\T_{+}(\H)\,|\,\Tr {G}\rho<+\infty\shs\}\,$ by the formula \footnote{We define the operator $A\rho B^*$ in such a way to avoid the notion of adjoint operator, since we make no assumptions about closability of the operators $A$ and $B$.}
\begin{equation}\label{ab-d}
  A\rho B^*\doteq\sum_i|\alpha_i\rangle\langle\beta_i|,\qquad |\alpha_i\rangle=A|\varphi_i\rangle,\;|\beta_i\rangle=B|\varphi_i\rangle,
\end{equation}
where $\rho=\sum_i |\varphi_i\rangle\langle\varphi_i|$ is any decomposition of $\rho\in\T^+_{G}$ into $1$-rank positive operators.}
\end{lemma}\medskip

\emph{Proof.} If $\rho=\sum_i |\varphi_i\rangle\langle\varphi_i|$ and $\{\psi_i\}$ is any  set of orthogonal unit vectors in a separable Hilbert space $\K$ then $|\eta\rangle=\sum_i |\varphi_i\rangle\otimes|\psi_i\rangle$ is a vector in $\,\D(\sqrt{G}\otimes I_{\K})$ such that $\rho=\Tr_{\K}|\eta\rangle\langle\eta|$. By Lemma \ref{p-l-1} the operators $A\otimes I_{\K}$
and $B\otimes I_{\K}$ have unique linear $\sqrt{G}\otimes I_{\K}$-bounded extensions to the set $\,\D(\sqrt{G}\otimes I_{\K})$ satisfying (\ref{s-prop}). Hence
\begin{equation}\label{lf}
\sum_i |A\varphi_i\rangle\langle B\varphi_i|=\Tr_{\K}|A\otimes I_{\K}\eta\rangle\langle B\otimes I_{\K}\eta|.
\end{equation}
So, by using the well known relation between different purifications of a given state \cite{H-SCI,Wilde} and Remark \ref{w-c}, it is easy to show that the r.h.s. of (\ref{lf}) does not depend on the representation $\rho=\sum_i |\varphi_i\rangle\langle\varphi_i|$. It follows that the r.h.s. of (\ref{lf})
correctly defines an affine function $\rho\mapsto A\rho B^*$ on the set $\T^+_{{G}}$. $\square$

Lemma \ref{p-l-2} implies, in particular, that $\rho\mapsto A\rho A^*$
is an affine function from $\T^+_{G}$ into $\T_+(\H)$ (well defined by formula (\ref{ab-d}) with $B=A$) for any $\sqrt{G}$-bounded operator $A$. Hence
the r.h.s. of (\ref{ec-on}) is well defined for any such operator.\smallskip

The following proposition shows that we may also define the operator \emph{E}-norm $\|A\|_E^G$ of any $\sqrt{G}$-bounded operator $A$ by formula (\ref{ec-on}).\smallskip

\begin{property}\label{v-new} \emph{Let $A$ be an arbitrary  $\sqrt{G}$-bounded operator and $E>0$.}\smallskip

\noindent A) \emph{The right hand sides of (\ref{ec-on-b}) and (\ref{ec-on}) coincide (provided that $A\rho A^*$ is defined by formula (\ref{ab-d}) with $B=A$).}

\noindent B) \emph{The suprema in definitions (\ref{ec-on-b}) and (\ref{ec-on}) can be taken, respectively, over all vectors in $\H_{\leq 1}$ satisfying the condition $\,\langle\varphi|G|\varphi\rangle\leq E$ and over all operators in $\T_{+,1}(\H)$ satisfying the condition $\,\Tr {G}\rho\leq E$.}
\end{property}\smallskip

\emph{Proof.} A) The concavity of the  function $E\mapsto \left[\|A\|^{G}_E\right]^2$ (Lemma \ref{new-l}A) and the proof of Proposition \ref{en-def}A
show that the supremum in the r.h.s. of (\ref{ec-on-f}) can be taken over all finite rank positive operators in $\T_{1,+}$ such that $\Tr G\rho\leq E$.

Let $\rho=\sum_{i=1}^{+\infty}|\varphi_i\rangle\langle\varphi_i|$ be an arbitrary state in $\S(\H)$ such that $\Tr G\rho\leq E$. Then $\rho_n=\sum_{i=1}^{n}|\varphi_i\rangle\langle\varphi_i|$
is a finite rank positive operator in $\T_{1,+}$ such that $\Tr G\rho_n\leq E$ for each $n$, and
$$
\lim_{n\rightarrow+\infty}\Tr A\rho_n A^*=\Tr A\rho A^*\leq+\infty.
$$
Thus, assertion A follows from equality (\ref{ec-on-f}) and the remark at the begin of the proof.

B) This assertion is proved by using concavity of the  function $E\mapsto \left[\|A\|^{G}_E\right]^2$ (Lemma \ref{new-l}A) and the proof of Proposition \ref{en-def}A. $\square$ \smallskip

By Proposition \ref{v-new}B for any  vector $\varphi$ in $\D(\sqrt{G})$ such that $\|\varphi\|\leq 1$ we have
\begin{equation}\label{a-phi-est}
 \|A\varphi\|\leq \|A\|^{G}_{E_{\varphi}}\leq K_{\varphi}\|A\|^{G}_E,\vspace{-5pt}
\end{equation}
where $\,E_{\varphi}=\|\sqrt{G}\varphi\|^2\,$ and $\,K_{\varphi}=\max\{1, \sqrt{E_{\varphi}/E}\}\,$. This implies the following
\smallskip
\begin{lemma}\label{vsl} \emph{Let $P_E$ be the spectral
projector of ${G}$ corresponding to the interval $[0,E]$. For any $\sqrt{G}$-bounded operator $A$ the operator
$AP_E$ is bounded and $\|AP_E\|\leq \|A\|^{G}_E$.}
\end{lemma}
\smallskip

\emph{Proof.} It follows from (\ref{a-phi-est}) that
$\|AP_E\varphi\|\leq \|A\|^{G}_E$ for any unit vector $\varphi$ in $\H$, since $\|\sqrt{G}P_E\varphi\|^2\leq E$ and $\|P_E\varphi\|\leq 1$. $\square$
\medskip

The following lemma shows that the set $\,\Pi_{\sqrt{G}}(A)$ is completely determined by the function $\,E\mapsto \left[\|A\|^{G}_E\right]^2$ and vice versa.\smallskip

\begin{lemma}\label{bl}
\emph{A pair $(a,b)$ belongs to the set $\,\Pi_{\sqrt{G}}(A)$ if and only if $\|A\|^{G}_{E}\leq \sqrt{a^2+b^2 E}$ for all $\,E>0$.}
\end{lemma}\smallskip

\emph{Proof.} If $\|A\|^{G}_{E}\leq \sqrt{a^2+b^2 E}$ then it follows from (\ref{a-phi-est}) that
$$
\|A\varphi\|\leq \|A\|^{G}_{\|\sqrt{G}\varphi\|^2}\leq \sqrt{a^2+b^2 \|\sqrt{G}\varphi\|^2}
$$
for any unit vector $\varphi$ in $\D(\sqrt{G})$. Hence $(a,b)\in\Pi_{\sqrt{G}}(A)$.
If $(a,b)\in\Pi_{\sqrt{G}}(A)$ then
$$
\sup\left\{\|A\varphi\|\,\left|\,  \varphi\in\D(\sqrt{G}), \|\varphi\|\leq1, \|\sqrt{G}\varphi\|^2\leq E\right\}\right.\leq\sqrt{a^2+b^2 E}
$$
for any $E>0$. So, definition (\ref{ec-on-b}) implies that $\|A\|^{G}_{E}\leq \sqrt{a^2+b^2 E}$. $\square$\smallskip

Denote by $\B_{\!G}(\H)$ the linear space of all $\sqrt{G}$-bounded operators equipped with the norm $\|\!\cdot\!\|^{G}_E$ defined by the equivalent expressions (\ref{ec-on-b}) and (\ref{ec-on}) (we identify operators coinciding on $\D(\sqrt{G})$). We will also consider the norm  $\sn\!\cdot\!\sn^{G}_E$ defined in (\ref{eq-norms-2}), which is commonly used on the space of $\sqrt{G}$-bounded operators.\smallskip

\begin{theorem}\label{comp-d} \emph{Let ${G}$ be a positive semidefinite unbounded densely defined operator on $\H$  satisfying condition (\ref{H-cond}).}\smallskip

\noindent A) \emph{$\B_{\!G}(\H)$ is a nonseparable Banach space. The norms  $\|\cdot\|^{G}_{E}$ and $\sn \cdot\sn^{G}_{E}$
satisfy the equivalence relations (\ref{E-n-eq}), (\ref{one}) and (\ref{E-n-eq+}) on $\,\B_{\!G}(\H)$ . For any $A\in\B_{\!G}(\H)$
and $E>0$ the following expressions hold}
$$
\sn A\sn^{G}_{E}=\sup_{t>0}\|A\|^{G}_{tE}/\sqrt{1+t},\qquad  \|A\|^{G}_{E}=\inf_{t>0}\sn A \sn^{G}_{tE}\sqrt{1+1/t}.
$$

\noindent B) \emph{If $A\in\B_G(\H)$ and $E>0$ then}\vspace{-5pt}
$$
\|A\|^{G}_{E}=\inf\left.\left\{\sqrt{a^2+b^2 E}\;\right|(a,b)\in\Pi_{\sqrt{G}}(A)\right\} \quad \textit{and} \quad b_{\sqrt{G}}(A)=\lim_{E\rightarrow+\infty}\|A\|^{G}_{E}/\sqrt{E}.
$$

\emph{The limit in the last formula can be replaced by the infimum over all $E>0$.} \smallskip\pagebreak

\noindent C) \emph{The completion of $\B(\H)$ w.r.t. any of the norms $\|\!\cdot\!\|^{G}_E$, $E>0$,
coincides with the closed subspace $\,\B^0_{\!G}(\H)$ of $\,\B_{\!G}(\H)$ consisting of all $\sqrt{G}$-infinitesimal operators, i.e. operators with
the $\sqrt{G}$-bound equal to $0$. An operator $A$ belongs to $\B^0_{\!G}(\H)$ if and only if}
\begin{equation}\label{s-cond}
 \|A\|^{G}_E=o\shs(\sqrt{E})\quad\textup{ as }\quad E\rightarrow+\infty.
\end{equation}

\emph{If $\,G$ is a discrete operator (Def.\ref{D-H}) then the Banach space $\B^0_{\!G}(\H)$ is separable.}\smallskip

\noindent D) \emph{Any ball in $\B(\H)$ is complete with respect to any of the norms $\|\!\cdot\!\|^{G}_E$, $E>0$. An operator $A$ belongs to $\B(\H)$ if and only if the function $E\mapsto \|A\|^{G}_E$ is bounded. In this case $\|A\|=\sup_{E>0}\|A\|^{G}_E=\lim_{E\rightarrow+\infty}\|A\|^{G}_E$.}\medskip

\noindent E) \emph{The $\sqrt{G}$-bound  is a continuous seminorm on $\B_G(\H)$. Quantitatively,
\begin{equation}\label{b-cb}
\left|\shs b_{\sqrt{G}}(A)-b_{\sqrt{G}}(B)\right|\leq b_{\sqrt{G}}(A-B)\leq \|A-B\|^{G}_{E}/\sqrt{E}
\end{equation}
for arbitrary $A,B$ in $\B_G(\H)$ and any $E>0$.}\smallskip

\noindent F) \emph{If $\,\K$ is a separable Hilbert space then  $\,\|A\otimes I_{\K}\|_E^{G\otimes I_{\K}}=\|A\|_E^{G}$  for any $A\in\B_G(\H)$}\footnote{$A\otimes I_{\K}$ denotes the operator mentioned in Lemma \ref{p-l-1}.}

\smallskip

\noindent G)  \emph{For arbitrary $\sqrt{G}$-bounded operators $A$ and $B$ and any operator $\rho$ in $\T_{+}(\H)$ such that $\Tr\rho\leq 1$ and $\,E_{\rho}\doteq\Tr G\rho<+\infty\,$ the following inequalities hold
$$
|\Tr A\rho B^*|\leq \|A\rho B^*\|_1 \leq \|A\|^{G}_{E_{\rho}}\|B\|^{G}_{E_{\rho}},
$$
where $A\rho B^*$ is the trace class operator defined in (\ref{ab-d})}.\smallskip



\noindent H) \emph{For any $A$ in $\B_{\!G}(\H)$ and $E>0$ the suprema in definitions (\ref{ec-on-b}) and (\ref{ec-on}) can be taken, respectively, over all unit vectors in $\H$ satisfying the condition
$\,\langle\varphi|G|\varphi\rangle=E$ and over all states in $\S(\H)$ satisfying the condition $\,\Tr {G}\rho=E$.}\smallskip

\noindent I) \emph{If the operator $G$ is  discrete (Def.\ref{D-H}) then for any $A\in\B^0_G(\H)$ and $E>0$ the suprema in (\ref{ec-on-b}) and (\ref{ec-on}) are attainable. Moreover, if $A$ is a unbounded operator in $\B^0_G(\H)$ or $\|A\|^{G}_{E}<\|A\|<+\infty$  then the suprema in (\ref{ec-on-b}) and (\ref{ec-on}) are attained, respectively, at unit vector $\varphi_0$ in $\H$ such that $\,\langle\varphi_0|G|\varphi_0\rangle=E$ and at a state
$\rho_0$ in $\S(\H)$ such that $\,\Tr G\rho_0=E$.}
\end{theorem}\medskip

\begin{remark}\label{comp-d-r+}
The expressions in Theorem 3A show that the functions $E\mapsto\|A\|^G_E$ and
$E\mapsto\sn A\sn^G_E$ are completely determined by each other for arbitrary $\sqrt{G}$-bounded operator $A$. So, if $\|A\|^G_E=\|B\|^G_E$ for all $E>0$ for some $\sqrt{G}$-bounded operators $A$ and $B$ then $\sn A\sn^G_E=\sn B\sn^G_E$ for all $E>0$ and vise versa.\footnote{This holds for the operators $a$ and $\sqrt{N}$ in the below example.} These expressions mean that the concave function $E\mapsto[\|A\|^G_E]^2$ and the function
$E\mapsto[\sn A\sn^G_E]^2$ are related by the transformations
$$
F[f](x)\ = \ \sup_{\, t > 0}\ \frac{f(xt)}{t+1}
\quad \textrm{and} \quad G[f](x)\ = \ \inf_{\, t > 0}\  f(xt)\, (1 \, + \, 1/t),
$$
defined on the set of nonnegative functions on $(0,+\infty)$. In \cite{P&Sh} it is shown that $G\circ F$ maps \emph{any} nonnegative function $f$ on  $(0,+\infty)$ into its concave hull and hence $G[F[f]]=f$ for any concave nonnegative function $f$. This shows that the second expression in Theorem 3A follows from the first one and the concavity of the function $E\mapsto[\|A\|^G_E]^2$.
\end{remark}\medskip

\begin{remark}\label{comp-d-r} The below proof of the density of $\B(\H)$ in $\B^0_{\!G}(\H)$ shows that
the  $\sqrt{G}$-infinitesimality criterion (\ref{s-cond}) is equivalent to the following one
\begin{equation}\label{s-cond+}
\lim_{n\rightarrow+\infty} \|A\bar{P}_n\|^{G}_E=0,
\end{equation}
where $\bar{P}_n$ is the spectral
projector of ${G}$ corresponding to the interval $\,(n,+\infty)$.
\end{remark}\smallskip

\emph{Proof of Theorem \ref{comp-d}.} A) By Lemma \ref{new-l}A inequalities (\ref{E-n-eq}), (\ref{one}) and (\ref{E-n-eq+}) for any $A$ in $\B_{\!G}(\H)$ are proved by the same arguments as for a bounded operator $A$.
By using definitions (\ref{ec-on-b}) and (\ref{eq-norms-2}) it is easy to show that
$$
\sn A\sn^{G}_{E}=\sup_{r\in(0,1)} \sqrt{r}\|A\|^{G}_{\frac{1-r}{r}E},\quad A\in\B_{\!G}(\H).
$$
The first expression in A follows from this one by the change of variables $t=(1-r)/r$.

The second expression in A is derived
from the first formula in part B proved below by noting that the infimum in that formula can be taken over all the pairs
$(\sn A\sn_E^G,\sn A\sn_E^G/\sqrt{E}), E>0$. This follows from density of the set
$$
\left.\left\{ \left(\sn A\sn_E^G+x, \sn A\sn_E^G/\sqrt{E}+y\right)\,\right|\, E>0,\, x,y\geq 0 \right\}
$$
in $\Pi_{\sqrt{G}}(A)$, which can be proved by noting that $\sn A\sn_E^G=\min\{ a\,|\,(a, a/\sqrt{E})\in \Pi_{\sqrt{G}}(A)\}$.

\smallskip

Denote by $\H_E^G$ the Hilbert space obtained by equipping the linear space $\D(\sqrt{G})$ with the inner product
$$
\langle\varphi|\psi\rangle_E^G=\langle\varphi|\psi\rangle+\langle\varphi|G|\psi\rangle/E.
$$

Since the norm $\sn A\sn_E^G$ of any $\sqrt{G}$-bounded operator $A$ is the operator norm of $A$ treated as a bounded operator from $\H_E^G$ into $\H$,
the linear space of all $\sqrt{G}$-bounded operators equipped with the norm $\sn \cdot\sn_E^G$ is a nonseparable Banach space \cite{R&S}.  Hence, the equivalence of the
norms $\sn \cdot\sn_E^G$ and $\|\cdot\|_E^G$ implies that $\B_{\!G}(\H)$ is a nonseparable Banach space. \smallskip

B) Since $E\mapsto [\|A\|^{G}_{E}]^2$ is a concave nonnegative function on $\mathbb{R}_+$, it coincides with the infimum of all linear functions $E\mapsto a^2+b^2 E$
such that $[\|A\|^{G}_{E}]^2\leq a^2+b^2 E$ for all $E>0$. The concavity of the function $E\mapsto [\|A\|^{G}_{E}]^2$ implies that the function $E\mapsto [\|A\|^{G}_{E}]^2/E$ is non-increasing. So, both formulae in part B
follow from Lemma \ref{bl}.\smallskip

C) The continuity and  the seminorm properties of the function $A\mapsto b_{\sqrt{G}}(A)$ stated in part E proved below show that
$\B^0_{\!G}(\H)=b^{-1}_{\sqrt{G}}(0)$  is a closed subspace of $\B_{\!G}(\H)$. The characterizing property (\ref{s-cond}) follows from the
second formula in part B.

To prove density of $\B(\H)$ in $\B^0_{\!G}(\H)$ it suffices, by Lemma \ref{vsl}, to show that for any $A\in\B^0_{\!G}(\H)$ the sequence $\{AP_n\}$, where $P_n$ is the spectral projector of $G$ corresponding to the interval $[0,n]$, converges to $A$ with respect to the norm $\|\!\cdot\!\|^{G}_E$.
For given $P_n$ let $\varphi$ be  any unit vector such that $\langle\varphi|G|\varphi\rangle\leq E$  and $x_n=\langle\varphi|\bar{P}_n|\varphi\rangle>0$, where  $\bar{P}_n=I_{\H}-P_n$. Let $|\varphi_n\rangle=x^{-1/2}_n \bar{P}_n|\varphi\rangle$. We have
$$
\|A\bar{P}_n \varphi\|^2=x_n\|A\varphi_n\|^2\leq x_n\left[\|A\|^{G}_{E/x_n}\right]^2\leq  (E/n)\left[\|A\|^{G}_{n}\right]^2.
$$
The first inequality follows from definition (\ref{ec-on-b}) of the operator \emph{E}-norm and the inequality
$\langle\varphi_n|G|\varphi_n\rangle\leq E/x_n$, the second one follows from concavity of the  function $E\mapsto \left[\|A\|^{G}_E\right]^2$, Lemma \ref{WL} and the inequality $x_n\leq E/n$ (which holds, since $\langle\varphi|G|\varphi\rangle\leq E$). The above estimate implies that
$$
\|A-AP_n\|^{G}_{E}\doteq \sup_{\varphi\in\mathcal{H}_1:
\langle\varphi|G|\varphi\rangle\leq E}\|A\bar{P}_n \varphi\|
\leq \sqrt{E/n}\|A\|^{G}_{n}.
$$
So, condition (\ref{s-cond}) guarantees that $\|A-AP_n\|^{G}_{E}$ tends to zero as $\,n\rightarrow+\infty$.\smallskip

The above arguments and Lemma \ref{vsl} imply that (\ref{s-cond}) is equivalent to (\ref{s-cond+}). \smallskip

If $G$ is a discrete operator then the separability of $\B^0_{\!G}(\H)$ follows from separability of $\B(\H)$ w.r.t. any of the operator \emph{E}-norms, which can be easily shown by using Proposition \ref{ec-on-p2}B
and separability of $\B(\H)$ w.r.t. the strong operator topology. \smallskip

D) We begin with the second assertion. The "only if" part of this  assertion and the expression $\|A\|=\sup_{E>0}\|A\|^{G}_E$ follow from Proposition \ref{ec-on-p1}. If  $\|A\|^{G}_E\leq M<+\infty$ for all $E>0$  then it follows from (\ref{a-phi-est}) that $\|A\varphi\|\leq M$ for any unit vector $\varphi$ in $\D(\sqrt{{G}})$. Since $\D(\sqrt{{G}})$ is dense in $\H$, this implies that $A\in\B(\H)$.\smallskip

To prove the first assertion assume that $\{A_n\}$ is a sequence in $\B(\H)$ converging to an operator $A_0\in\B^0_{\!G}(\H)$ such that $\|A_n\|\leq M<+\infty$ for all $n$. Since $\,\|A_n\|^{G}_E\leq\|A_n\|\leq M$ for all $n$ and $E>0$ and the right hand side of the inequality
$$
\left|\|A_n\|^{G}_E-\|A_0\|^{G}_E\right|\leq\|A_n-A_0\|^{G}_E
$$
tends to zero as $n\rightarrow+\infty$ for any $E>0$, it is easy to see that $\|A_0\|^{G}_E\leq M$ for all $E$. Thus, $\|A_0\|\leq M$ by the assertion proved before. \smallskip

E) The seminorm properites of $b_{\sqrt{G}}(\cdot)$ follow from the second formula in part B of the theorem. So, since the function $E\mapsto [\|A\|^{G}_{E}]^2/E$ is non-increasing for any given $A\in\B_{\!G}(\H)$, the inequality (\ref{b-cb}) follows from the triangle inequality for $b_{\sqrt{G}}(\cdot)$.\smallskip

F) This assertion follows from Lemma \ref{p-l-1} and the first formula in part B of the theorem. \smallskip

G) Let $\rho=\sum_i |\varphi_i\rangle\langle\varphi_i|$ be a decomposition into $1$-rank positive operators and $\{\psi_i\}$  a set of orthogonal unit vectors in a separable Hilbert space $\K$ then $|\eta\rangle=\sum_i |\varphi_i\rangle\otimes|\psi_i\rangle$ is a vector in $\,\D(\sqrt{G}\otimes I_{\K})$ such that $\rho=\Tr_{\K}|\eta\rangle\langle\eta|$. By Lemma \ref{p-l-1} the operators $A\otimes I_{\K}$
and $B\otimes I_{\K}$ have unique $\sqrt{G}\otimes I_{\K}$-bounded linear extensions to the set $\,\D(\sqrt{G}\otimes I_{\K})$ satisfying (\ref{s-prop}). By the monotonicity of the trace norm we have
$$
\|A\rho B^*\|_1\leq\||A\otimes I_{\K}\eta\rangle\langle B\otimes I_{\K}\eta|\|_1\leq\|A\otimes I_{\K}\eta\|\| B\otimes I_{\K}\eta\|\leq\|A\otimes I_{\K}\|^{{G}\otimes I_{\K}}_{E_{\rho}}\|B\otimes I_{\K}\|^{{G}\otimes I_{\K}}_{E_{\rho}}.
$$
By part F of the theorem  the r.h.s. of this inequality is equal to $\|A\|^{{G}}_{E_{\rho}}\|B\|^{{G}}_{E_{\rho}}$. \smallskip

H) If $A$ is a bounded operator then the possibility to take the suprema in (\ref{ec-on-b}) and (\ref{ec-on}) over all unit vectors in $\H$ satisfying the condition
$\,\langle\varphi|G|\varphi\rangle=E$ and over all states in $\S(\H)$ satisfying the condition $\,\Tr {G}\rho=E$ correspondingly follows Proposition \ref{en-def}B. If $A$ is a unbounded operator then this possibility can be easily shown by noting that the function $E\mapsto\|A\|^{G}_E$ is strictly increasing on $\mathbb{R}_+$ (since it is concave on $\mathbb{R}_+$ and tends to $+\infty$ as $E\rightarrow+\infty$). \smallskip

I) If $A\in\B_{\!G}^0(\H)$ then $\rho\mapsto \Tr A\rho A^*$ is an affine continuous function on the set $\mathfrak{C}_{{G},E}\doteq \{\rho\in\T_{+}(\H)\,|\,\Tr\rho\leq 1, \Tr {G}\rho\leq E\shs\}$ for any $E>0$ by Corollary \ref{vbl-c}
below (proved independently). So, both assertions are proved by repeating the arguments from the proof of Proposition \ref{en-def}C. $\square$
\smallskip

\textbf{Example: the operators associated with the Heisenberg
Commutation Relation}

Let $\H = L_2(\mathbb{R})$  and $S(\mathbb{R})$ be the set of infinitely differentiable rapidly decreasing functions
with all the derivatives tending to zero quicker than any degree of $|x|$ when $|x|\rightarrow+\infty$.
Consider the operators $q$ and $p$ defined on the set $S(\mathbb{R})$ by setting
$$
(q\varphi)(x) = x\varphi(x)\quad\textrm{and}\quad  (p\shs\varphi)(x) = \frac{1}{i}\frac{d}{dx}\varphi(x).
$$
These operators are essentially self-adjoint. They represent (sharp) real
observables of position and momentum of a quantum particle in the system of units where Planck's
constant $\hbar$ equals to $1$ \cite[Ch.12]{H-SCI}. On the domain $S(\mathbb{R})$ these operators  satisfy the Heisenberg commutation relation
\begin{equation}\label{H-cr}
[q, p\shs] = i I_{\H}.
\end{equation}

For given $\omega>0$ consider the operators
\begin{equation}\label{a-oper-d}
a=(\omega q+ip)/\sqrt{2\omega}\quad \textrm{and} \quad a^{\dagger}=(\omega q-ip)/\sqrt{2\omega}
\end{equation}
defined on $S(\mathbb{R})$. Via these operators the commutation relation (\ref{H-cr}) can be rewritten as
$
[a,a^{\dagger}] = I_{\H}.
$
The operator
\begin{equation}\label{a-N}
N=a^{\dagger}a=aa^{\dagger}-I_{\H}
\end{equation}
is positive and essentially self-adjoint. It represents (sharp) real
observable of the number of quanta of the harmonic oscillator with the frequency $\omega$.
The selfadjoint extension of $N$ has the form (\ref{H-rep}) with $E_n=n$ and the basic $\{\tau_n\}$ of eigenvectors of $N$ which
can be described as follows
$$
\tau_0(x)=\sqrt[4]{\frac{\omega}{\pi}}\exp\left[-\frac{\omega x^2}{2}\right],\quad |\tau_n\rangle=\frac{1}{\sqrt{n!}}\,[a^{\dagger}]^n|\tau_{0}\rangle,\; n\geq 1.
$$
So, $N$ is a positive unbounded discrete (Def.\ref{D-H}) operator satisfying condition (\ref{H-cond}).

The operators $a$ and $a^{\dagger}=a^*$ are called \emph{annihilation} and \emph{creation} operators correspondingly, since
\begin{equation}\label{a-oper}
a|\tau_0\rangle=0,\quad a|\tau_n\rangle=\sqrt{n}|\tau_{n-1}\rangle\quad\textrm{and}\quad a^{\dagger}|\tau_n\rangle=\sqrt{n+1}|\tau_{n+1}\rangle.
\end{equation}
So, the operators $a$ and $a^{\dagger}$ are correctly extended to the set
$$
\D(\sqrt{N})=\left\{\varphi\in\H\,\left|\,\sum_{n=0}^{\infty}n |\langle\varphi|\tau_n\rangle|^2<+\infty\right.\right\}
$$
By using relations (\ref{a-oper-d}) the operators $p$ and $q$ are also extended to the set $\D(\sqrt{N})$.

We will estimate the operator \emph{E}-norm of the operators $q$, $p$, $a$ and $a^{\dagger}$
induced by the operator $N$ (which up to the constant summand coincides with the Hamiltonian of a quantum oscillator).
By using (\ref{a-N}) it is easy to show that $\|a\|_E^N=\|\sqrt{N}\|_E^N=\sqrt{E}$ and $\|a^{\dagger}\|_E^N=\sqrt{E+1}$ for any $E>0$.
For the operators $\,q=(a^{\dagger}+a)/\sqrt{2\omega}\,$ and $\,p=i\sqrt{\omega/2}(a^{\dagger}-a)\,$
one can obtain the following estimates
\begin{equation}\label{QP-en}
\sqrt{\frac{2E+1/2}{\omega}}< \|q\|^{N}_E\leq\sqrt{\frac{2E+1}{\omega}},\quad \sqrt{(2E+1/2)\omega}< \|p\|^N_E\leq\sqrt{(2E+1)\omega}
\end{equation}
(the \emph{E}-norms of $q$ and $p$ depend on $\omega$, since the operator $N$ depends on $\omega$). The right inequalities in (\ref{QP-en}) directly follow from
the triangle inequality and the above expressions for $\|a\|_E^N$ and $\|a^{\dagger}\|_E^N$. To prove the left inequalities in (\ref{QP-en}) it suffices to show that
$$
\sup_{\substack{\|\varphi\|=1,\shs\langle\varphi|N|\varphi\rangle\leq E}}\|(a^{\dagger}\pm a)\varphi\|>\sqrt{4E+1}.
$$
This can be easily done by using the unit vectors $|\varphi_{\pm}\rangle= \sqrt{1-r}\sum_{n=0}^{+\infty} c^{\pm}_n r^{n/2}|\tau_n\rangle$, where $r=E/(E+1)$, $c_n^{-}=e^{i\pi n/2}$ and $c^+_n=1$ for all $n$.\smallskip

By using the first expression in Theorem \ref{comp-d}A and the above estimates of the norms $\|a\|_E^N$, $\|a^{\dagger}\|_E^N$, $\|p\|_E^N$ and $\|q\|_E^N$
we obtain $\sn a\sn_E^N=\sqrt{E}$, $\sn a^{\dagger}\sn_E^N=\max\{1,\sqrt{E}\}$,
\begin{equation*}
\sqrt{\frac{l(E)}{\omega}}< \sn q\sn^{N}_E\leq\sqrt{\frac{u(E)}{\omega}}\quad \textrm{and} \quad \sqrt{l(E)\omega}< \sn p\sn^N_E\leq\sqrt{u(E)\omega},
\end{equation*}
where $l(E)=\max\{1/2, 2E\}$ and $u(E)=\max\{1, 2E\}$.\smallskip

The second formula in Theorem \ref{comp-d}B and the above estimates of the norms $\|a\|_E^N,$ $\|a^{\dagger}\|_E^N$, $\|p\|_E^N$ and $\|q\|_E^N$ imply that
$$
b_{\sqrt{N}}(a)=b_{\sqrt{N}}(a^{\dagger})=1,\qquad b_{\sqrt{N}}(q)=\sqrt{2/\omega}\quad \textrm{and}\quad b_{\sqrt{N}}(p)=\sqrt{2\omega}.
$$

So, the operators $q$, $p$, $a$ and  $a^{\dagger}$ belong to the Banach space $\B_{\!N}(\H)$ but not lie in the completion
$\B^0_{\!N}(\H)$ of $\B(\H)$ w.r.t. the norm $\|\!\cdot\!\|^N_E$.\smallskip

For any $t<1$ let $a_t$ and $a_t^{\dagger}$ be the operators defined on the set $\D(\sqrt{N})$ by settings
\begin{equation}\label{a-oper}
a_t|\tau_0\rangle=0,\quad a_t|\tau_n\rangle=n^{t/2}|\tau_{n-1}\rangle\quad\textrm{and}\quad a_t^{\dagger}|\tau_n\rangle=(n+1)^{t/2}|\tau_{n+1}\rangle.
\end{equation}
It is easy to show that
\begin{equation}\label{a-soc}
\lim_{t\rightarrow 1-0} a_t|\varphi\rangle=a|\varphi\rangle\quad \textrm{and} \quad
\lim_{t\rightarrow 1-0} a_t^{\dagger}|\varphi\rangle=a^{\dagger}|\varphi\rangle\quad \textrm{for any}\;\; \varphi\in\D(\sqrt{N}).
\end{equation}

Since $a_t^{\dagger}a_t=N^t$ and $a_ta_t^{\dagger}=(N+I_{\H})^t$, by using concavity of the function $x\mapsto x^{t}$, we obtain
$$
\|a_t\|_E^N\leq \sqrt{\sup_{\Tr N\rho\leq E } [\Tr N \rho\shs]^t}=E^{t/2},\quad  \|a^{\dagger}_t\|_E^N\leq \sqrt{\sup_{\Tr N\rho\leq E } [\Tr (N+I_{\H}) \rho\shs]^t}=(E+1)^{t/2}.
$$

So, the operators $a_t$ and $a_t^{\dagger}$ belong to the space $\B^0_{\!N}(\H)$ for all $t<1$ (since they satisfy condition (\ref{s-cond})), while the "limit" operators $a$ and $a^{\dagger}$ lie in  $\B_{\!N}(\H)\setminus\B^0_{\!N}(\H)$. So, $a_t$ and $a_t^{\dagger}$ do not tend to $a$ and $a^{\dagger}$ as $t\rightarrow 1$ w.r.t. the norm $\|\!\cdot\!\|_E^N$ in spite of the strong operator convergence (\ref{a-soc}).
\smallskip

\begin{remark}\label{thelast}
It follows from (\ref{a-phi-est}) that
\begin{equation}\label{vlast-eq}
\|\!\cdot\!\|^{G}_E\,\textrm{-}\lim_{n\rightarrow\infty}A_n=A_0\quad\Rightarrow\quad\lim_{n\rightarrow\infty}A_n|\varphi\rangle=A_0|\varphi\rangle \quad\forall \varphi\in\D(\sqrt{G})
\end{equation}
for a sequence $\{A_n\}\subset\B_{\!G}(\H)$. The above example shows that the converse implication is not valid even in the case of discrete
operator $G$ (in this case $"\Leftrightarrow"$ holds in (\ref{vlast-eq}) for any bounded sequence $\{A_n\}\subset\B(\H)$ by Proposition 2B).
\end{remark}\medskip

In the last part of this section we consider properties of the Banach space $\B^0_{\!G}(\H)$.\smallskip

\begin{property}\label{vbl} \emph{If $A\in\B^0_{\!G}(\H)$ then the extension of $A\otimes I_{\K}$ to the set $\D(\sqrt{G}\otimes I_{\K})$ mentioned in Lemma \ref{p-l-1} is uniformly  continuous on the set
\begin{equation}\label{s-3}
\V_{E}\doteq \{\shs\eta\in\D(\sqrt{G}\otimes I_{\K})\,|\, \|\sqrt{G}\otimes I_{\K}\eta\|^2\leq E\shs\}
\end{equation}
for any $E>0$. Quantitatively,
\begin{equation}\label{v-CB}
\|A\otimes I_{\K}(\eta-\theta)\|\leq \varepsilon\|A\|^{G}_{4E/\varepsilon^2}=o(1)\quad \textit{as} \;\; \varepsilon\to 0^+
\end{equation}
for any vectors $\,\eta$ and $\,\theta$ in $\V_{E}$ such that $\|\eta-\theta\|\leq\varepsilon$.}\smallskip

\emph{If $\,A\in\B_{\!G}(\H)\setminus\B^0_{\!G}(\H)$ then the operator $A\otimes I_{\K}$ is not continuous  on the set $\,\V_{E}$ for any $E>0$.}\smallskip
\end{property}\medskip

\emph{Proof.} By Theorem \ref{comp-d}F for any unit  vector $\eta$ in $\D(\sqrt{G}\otimes I_{\K})$ we have
\begin{equation}\label{gen-ineq}
  \|A\otimes I_{\K}\, \eta\|\leq \|A\otimes I_{\K}\|^{{G}\otimes I_{\K}}_{E_{\eta}}=\|A\|^{G}_{E_{\eta}},\quad \textrm{where}\quad E_{\eta}=\|\sqrt{G}\otimes I_{\K} \eta \|^2.
\end{equation}

Assume that $\eta$ and $\theta$ are vectors in $\V_E$ such that $\|\eta-\theta\|\leq\varepsilon$. Since $\|\sqrt{G}\otimes I_{\K}\eta\|^2\leq E$ and $\|\sqrt{G}\otimes I_{\K}\theta\|^2\leq E$ we have $\|\sqrt{G}\otimes I_{\K}(\eta-\theta)\|^2\leq 4E$. So, by using (\ref{gen-ineq}), the concavity of the function $E\mapsto \left[\|A\|_E^{G}\right]^2$ and Lemma \ref{WL} we obtain
\begin{equation*}
\!\|A\otimes I_{\K} (\eta-\theta)\|=\|\eta-\theta\|\left\|A\otimes I_{\K}\, \frac{\eta-\theta}{\|\eta-\theta\|}\right\|\leq\|\eta-\theta\|\|A\|^{G}_{4E/\|\eta-\theta\|^2}\leq\varepsilon\|A\|^{G}_{4E/\varepsilon^2}.
\end{equation*}
By condition (\ref{s-cond}) the r.h.s. of this inequality tends to zero as $\varepsilon\rightarrow0^+$. Thus, the function $\eta\mapsto A\otimes I_{\K}|\eta\rangle$ is uniformly continuous on $\V_E$. \smallskip

The last assertion of the proposition follows from the proof of the last assertion of Corollary \ref{vbl-c} below, since
$\Tr_{\K}|A\otimes I_{\K}\eta\rangle\langle A\otimes I_{\K}\eta|=\Tr A\rho_{\eta}A^*$, where $\rho_{\eta}=\Tr_{\K}|\eta\rangle\langle\eta|$,
for any vector $\eta$ in $\D(\sqrt{G}\otimes I_{\K})$. $\square$ \medskip


\begin{corollary}\label{vbl-c} \emph{For any operators $A$ and $B$ in $\B_{\!G}^0(\H)$
the function $\shs\rho\mapsto A\rho B^*$ from $\T^+_{G}$ into $\T(\H)$ (defined by formula (\ref{ab-d})) is uniformly  continuous on the set
$\mathfrak{C}_{{G},E}\doteq \{\rho\in\T_{+}(\H)\,|\,\Tr\rho\leq 1,  \Tr {G}\rho\leq E\shs\}$ for any $E>0$. Quantitatively,
\begin{equation}\label{ab-cb}
\|A\rho B^*- A\sigma B^*\|_1\leq \sqrt{\varepsilon}\left(\|A\|_E^G \|B\|^{G}_{4E/\varepsilon}+\|B\|_E^G\|A\|^{G}_{4E/\varepsilon}\right)=o(1)\quad \textit{as} \;\; \varepsilon\to 0^+
\end{equation}
for any operators $\,\rho$ and $\sigma$ in $\,\C_{{G},E}$ such that $\|\rho-\sigma\|_1\leq\varepsilon$.}

\emph{If $\,A\in\B_{\!G}(\H)\setminus\B^0_{\!G}(\H)$ then  the function  $\,\rho\mapsto A\rho A^*$ is not continuous on the set $\,\mathfrak{C}_{{G},E}$  for any $E>0$.}
\medskip
\end{corollary}\medskip

\begin{remark}\label{comp-d-r++}
Corollary \ref{vbl-c} shows that the operators $A$ in $\B^0_{\!G}(\H)$ are characterized by continuity of the function $\rho\mapsto A\rho A^*$ on the set $\mathfrak{C}_{{G},E}$ for any given $E>0$.
\end{remark}\smallskip

\emph{Proof.} Let $\rho$ and $\sigma$ be operators in $\mathfrak{C}_{{G},E}$ such that $\|\rho-\sigma\|_1\leq\varepsilon$. If $\K\cong\H$ then
one can find  vectors $\eta$ and $\theta$ in the set $\V_{E}$ (defined in (\ref{s-3})) such that $\rho=\Tr_{\K} |\eta\rangle\langle\eta|$, $\sigma=\Tr_{\K} |\theta\rangle\langle\theta|$ and $\|\eta-\theta\|\leq \sqrt{\varepsilon}$ \cite{H-SCI}. By Lemma \ref{p-l-1} the operators $A\otimes I_{\K}$
and $B\otimes I_{\K}$ have unique linear $\sqrt{G}\otimes I_{\K}$-bounded extensions to the set $\,\D(\sqrt{G}\otimes I_{\K})$ satisfying (\ref{s-prop}). By using the monotonicity of the trace norm, the inequality $$\||\alpha\rangle\langle\beta|-|\varphi\rangle\langle\psi|\|_1\leq \|\alpha\|\|\beta-\psi\|+\|\psi\|\|\alpha-\varphi\|,$$  where $|\alpha\rangle=A\otimes I_{\K}|\eta\rangle$, $|\beta\rangle=B\otimes I_{\K}|\eta\rangle$, $|\varphi\rangle=A\otimes I_{\K}|\theta\rangle$, $|\psi\rangle=B\otimes I_{\K}|\theta\rangle$, and continuity bound (\ref{v-CB}) we obtain
$$
\|A\rho B^*-A\sigma B^*\|_1\leq \sqrt{\varepsilon}\|B\|^{G}_{4E/\varepsilon}\|A\otimes I_{\K}\shs\eta\|+\sqrt{\varepsilon}\|A\|^{G}_{4E/\varepsilon}\|B\otimes I_{\K}\shs\theta\|.
$$
By inequality (\ref{gen-ineq}) this implies (\ref{ab-cb}).

The r.h.s. of (\ref{ab-cb}) tends to zero as $\,\varepsilon\rightarrow 0^+$, since $A$ and $B$ satisfy condition (\ref{s-cond}).

If $A\in\B_{\!G}(\H)\setminus\B^0_{\!G}(\H)$ then, by Remark \ref{comp-d-r}, the sequence
$\|A\bar{P}_n\|^{G}_E$, where $\bar{P}_n$  is the spectral
projector of ${G}$ corresponding to the interval $\,(n,+\infty)$, does not tend to zero. Hence there is a sequence $\{\rho_n\}$ of states in $\mathfrak{C}_{{G},E}$ such that the sequence $\{\Tr A\bar{P}_n\rho_n\bar{P}_n A^*\}$ does not tend to zero. Since the condition $\Tr {G}\rho_n\leq E$ implies $\Tr \bar{P}_n\rho_n\leq E/n$, the sequence $\{\bar{P}_n\rho_n\bar{P}_n\}\subset\mathfrak{C}_{{G},E}$ tends to zero. This shows discontinuity of the function $\rho\mapsto A\rho A^*$ on the set $\mathfrak{C}_{{G},E}$.
$\square$
\pagebreak

Since $\B(\H)$ is dense in $\B_{\!G}^0(\H)$ by Theorem \ref{comp-d}C, Proposition \ref{bp-en-0}E implies the following\smallskip

\begin{property}\label{bp-en-0+} \emph{Let ${G}$ be a positive  densely defined operator on $\H$ satisfying condition (\ref{H-cond}) and $E>0$.}
\emph{Any 2-positive linear map $\Phi:\B(\H)\rightarrow\B(\H)$ such that $\,\Phi(I_{\H})\leq I_{\H}\,$ having the predual map $\Phi_*:\T(\H)\rightarrow\T(\H)$
with finite\footnote{By concavity of the function $\,E\mapsto Y_{\Phi}(E)\,$ its finiteness for some $E>0$ implies its finiteness for all $E>0$.}
\begin{equation}\label{E-phi}
Y_{\Phi}(E)\doteq \sup\left\{\shs\Tr {G}\Phi_*(\rho)\,|\,\rho\in\mathfrak{S}(\mathcal{H}),\Tr {G}\rho\leq E\,\right\}
\end{equation}
is uniquely extended to  the bounded linear operator $\,\Phi_{\!G}:\B^0_{\!G}(\H)\rightarrow\B_{\!G}^0(\H)$ such that}
\begin{equation}\label{ch-ext}
\|\Phi_{\!G}(A)\|^{G}_E\leq\sqrt{\|\Phi(I_{\H})\|}\,\|A\|^{G}_{Y_{\Phi}(E)}\leq \sqrt{\|\Phi(I_{\H})\|K_{\Phi}}\,\|A\|^{G}_{E},
\end{equation}
\emph{where }$K_{\Phi}=\max\{1, Y_{\Phi}(E)/E\}$.
\end{property}\smallskip

The assertion of Proposition \ref{bp-en-0+} can be strengthened substantially by assuming complete positivity of $\Phi$. The corresponding result is  considered in \cite{EPM}.\medskip

Different applications of the operator \emph{E}-norms  are presented in \cite{SPM,QDS,EPM}. In \cite{SPM} the version of the Kretschmann-Schlingemann-Werner theorem for unbounded completely positive linear maps is obtained by using the results from Section 5.
\bigskip\bigskip

I am grateful to A.S.Holevo, G.G.Amosov, A.V.Bulinsky and M.M.Wilde for discussion and useful remarks.
I am also grateful to V.Zh.Sakbaev for consultation concerning unbounded operators and to T.V.Shulman for the help and useful discussion.

Special thanks to S.Weis for the help in proving the coincidence of definitions (\ref{ec-on-b}) and (\ref{ec-on}).

\end{document}